\pdfoutput=1

\documentclass[final,3p,12pt,authoryear,hidelinks,times]{elsarticle}
\usepackage{placeins} %

\usepackage[utf8]{inputenc}

\usepackage[font=normalsize,labelfont=bf]{caption}

\usepackage{hyperref}
\hypersetup{colorlinks=false}
\usepackage[displaymath]{lineno}

\journal{Computers \& Structures}

\biboptions{authoryear}

\usepackage{amssymb,amsmath}
\usepackage{bm}

\usepackage{textcomp}
\usepackage{gensymb}

\usepackage{subcaption}

\usepackage{longtable,booktabs,tabularx,tabulary}
\newcolumntype{L}[1]{>{\hsize=#1\hsize\raggedright\arraybackslash}X}%
\newcolumntype{R}[1]{>{\hsize=#1\hsize\raggedleft\arraybackslash}X}%
\newcolumntype{C}[2]{>{\hsize=#1\hsize\columncolor{#2}\centering\arraybackslash}X}%
\newcolumntype{Y}{>{\centering\arraybackslash}X}

\newcommand*\patchAmsMathEnvironmentForLineno[1]{%
  \expandafter\let\csname old#1\expandafter\endcsname\csname #1\endcsname
  \expandafter\let\csname oldend#1\expandafter\endcsname\csname end#1\endcsname
  \renewenvironment{#1}%
     {\linenomath\csname old#1\endcsname}%
     {\csname oldend#1\endcsname\endlinenomath}}%
\newcommand*\patchBothAmsMathEnvironmentsForLineno[1]{%
  \patchAmsMathEnvironmentForLineno{#1}%
  \patchAmsMathEnvironmentForLineno{#1*}}%
\AtBeginDocument{%
\patchBothAmsMathEnvironmentsForLineno{equation}%
\patchBothAmsMathEnvironmentsForLineno{align}%
\patchBothAmsMathEnvironmentsForLineno{flalign}%
\patchBothAmsMathEnvironmentsForLineno{alignat}%
\patchBothAmsMathEnvironmentsForLineno{gather}%
\patchBothAmsMathEnvironmentsForLineno{multline}%
}

\begin{document}

\begin{frontmatter}

\title{Recovery by discretization corrected\\ particle strength exchange (DC PSE) operators}

\address[uwa]{Intelligent Systems for Medicine Laboratory,
The University of Western Australia,\\
35 Stirling Highway, Perth, Western Australia}

\address[harvard]{Harvard Medical School, Boston, Massachusetts,
USA}

\author[uwa]{B.F. Zwick\corref{mycorrespondingauthor}}
\ead{benjamin.zwick@uwa.edu.au}

\author[uwa]{G.C. Bourantas}

\author[uwa]{F. Alkhatib}

\author[uwa]{A. Wittek} 

\author[uwa,harvard]{K. Miller}

\cortext[mycorrespondingauthor]{Corresponding author}

\begin{abstract}
  A new recovery technique based on discretization corrected
  particle strength exchange (DC PSE) operators
  is developed in this paper.
  DC PSE is a collocation method that can be used
  to compute derivatives directly at nodal points,
  instead of by projection from Gauss points
  as is done in many finite element-based recovery techniques.
  The proposed method is truly meshless
  and does not require patches of elements to be defined,
  which makes it generally applicable to
  point clouds and arbitrary element topologies.
  Numerical examples show that the proposed method
  is accurate and robust.
\end{abstract}

\begin{keyword}
  meshless and meshfree methods\sep
  finite element method\sep
  linear elasticity\sep
  solid mechanics\sep
  numerical differentiation\sep
  patch recovery
\end{keyword}

\end{frontmatter}

\section{Introduction}

In computational science and engineering
there are many problems in which the gradient of the solution,
rather than the solution itself,
is of primary interest.
Examples include strain and stress in structural analysis
and wall shear stress in fluid dynamics.
The goal of this paper is to develop a meshfree method
for the accurate and robust recovery of gradients,
and other derived quantities such as strain and stress,
from discrete scalar or vector fields with values defined at nodal points.
These nodal values are typically obtained numerically,
using for example the finite element (FE) method.
However, the nodal values may also be obtained
analytically or from measurements, %
an example of the latter being the computation of strain
from image data
\citep{choi_etal_1991_image,lopez-linares_etal_2019_imagebased}. %

The finite element method (FEM)
has become the dominant method for solving partial differential equations (PDEs)
in problems with complex geometries.
However, the calculation of accurate stresses in the displacement-based FE method
is a notoriously difficult problem
\citep{oden_brauchli_1971_calculation}.
In the FE method, the function (e.g.\ displacement)
is best sampled at the nodal points
but the best accuracy for the gradient
(and derived quantities such as strain and stress)
is obtained at the superconvergent Gauss points
corresponding to the order of polynomial used for the solution
\citep{zienkiewicz_zhu_1995_superconvergence,zienkiewicz_etal_2013_finite_solid}.
Recovery methods are often used to smooth the discontinuous derivatives
and derived quantities (strain and stress)
resulting from piecewise continuous approximations of the solution variables
in finite element computations
\citep{zienkiewicz_etal_2013_finite_solid}.
Numerous recovery methods have been proposed
to obtain continuous gradient, strain and stress fields
from the discontinuous fields obtained using FEM
\citep{oden_brauchli_1971_calculation,
kelly_etal_1983_posteriori,
zienkiewicz_zhu_1987_simple,
zienkiewicz_zhu_1992_superconvergent_part1_recovery,
boroomand_zienkiewicz_1997_recovery,
zhang_naga_2005_new}.
An important application of recovery methods
is their use in \emph{a posteriori} error estimators
that are based on measuring the difference between
the direct and post-processed approximations to the gradient
\citep{zienkiewicz_zhu_1987_simple,ainsworth_oden_2000_posteriori}.

Meshless and meshfree methods
\citep{chen_belytschko_2015_meshless,fasshauer_2007_meshfree}
have been proposed as an alternative
to global and local recovery methods based on finite elements.
Compared with finite element based methods
which often require patches of neighboring elements to be defined,
meshfree methods do not rely on element connectivity
and can therefore be applied easily to almost arbitrary point clouds
(some limitations on the distribution of points remain
based on invertibility of matrices and conditioning of linear systems that must be solved
to construct the shape functions).
The problem is typically posed as an approximation or curve fitting problem
whereby the values are known at one set of points
(typically the quadrature points)
and these are to be interpolated as accurately as possible
to a different set of points
(typically the nodal points).
This approach can also be applied to the recovery of stress
in meshfree solution methods such as the element-free Galerkin method
\citep{lee_zhou_2004_error}.
\citet{ahmed_etal_2018_interpolation}
showed that a meshfree approach
based on moving least squares approximation
can be used to improve the performance of error recovery techniques
as compared to the superconvergent patch recovery (SPR) method
\citep{zienkiewicz_zhu_1992_superconvergent_part1_recovery}.
\citet{ahmed_2020_comparative}
compared the element-free Galerkin and radial point interpolation methods 
that use the moving least squares (MLS) interpolation (radial weighting function) and the radial point interpolation (polynomial and radial basis function), respectively,
against the mesh-dependent least squares interpolation error recovery (SPR)
\citep{zienkiewicz_zhu_1992_superconvergent_part1_recovery}
and found that the quality of meshless and mesh-based RPIM and MLS based error estimation
is superior to mesh-dependent SPR.

In this paper,
rather than evaluating the gradient, strain or stress at the integration points
and interpolating these values to the nodal points,
we aim to compute the gradient, strain and stress
directly at the nodal points given the nodal data.
To this end, we develop
an accurate, efficient and generally applicable meshfree approach
for recovery of gradients given the solution at nodal points
based on discretization corrected particle strength exchange (DC PSE) operators
\citep{schrader_etal_2010_discretization}.
DC PSE uses Taylor expansions on irregularly distributed point clouds to compute derivatives.
We compute the derivatives of the unknown field function
as a postprocessing procedure
after its solution has been obtained
and use these derivatives to compute derived quantities
such as strain and stress.

The remainder of the paper is organized as follows.
In section~\ref{sec:methods},
we state the governing equations of linear elasticity
and describe their solution using the FEM.
We then present the DC PSE operators
as an alternative to the FEM
for computing derivatives at nodal points,
and explain how we apply DC PSE as a post-processing tool
to obtain strain and stress from a given displacement field.
In section~\ref{sec:examples},
we demonstrate the accuracy of the proposed recovery method
by applying it to several benchmark problems,
as well as a more realistic practical example
involving the calculation of stress within the wall
of an abdominal aortic aneurysm (AAA).
Finally, in section~\ref{sec:discussion}
we discuss the pros and cons of the proposed approach
and make some suggestions for further research in this area.

\FloatBarrier
\section{Methods}
\label{sec:methods}

\subsection{Governing equations of linearized elasticity}

In this study, we consider elliptic PDEs of the following form.
Consider the region
$\Omega \subset \mathbb{R}^d$,
where $d = 1, 2$ (plane strain) or $3$.
The boundary \(\partial\Omega\) of the domain
is decomposed into two disjoint parts such that
\(\overline{\Gamma_D \cup \Gamma_N} = \partial\Omega\) and
\(\Gamma_D \cap \Gamma_N = \emptyset\).
The unit outward normal vector on \(\partial\Omega\)
is denoted by \(\bm{n}\).
Known displacements and tractions,
$\bm{u}_D$ and $\bm{g}$,
are prescribed on \(\Gamma_D\) and \(\Gamma_N\), respectively.
Body forces $\bm{f}$ are prescribed in \(\Omega\).
The strong form of the boundary value problem of elastostatics
can then be written as follows.
Given
\(\bm{f}:\Omega \to \mathbb{R}^d\),
\(\bm{g}:\Gamma_N \to \mathbb{R}^d\) and
\(\bm{u}_D:\Gamma_D \to \mathbb{R}^d\);
find
\(\bm{u}:\Omega \to \mathbb{R}^d\)
such that
\begin{alignat}{2}
  \label{eq:elasticity}
  -\mathrm{div} \, \bm{\sigma}(\bm{u}) &= \bm{f}   &\quad &\text{ in }\Omega, \\
  \bm{\sigma} (\bm{u}) \cdot \bm{n}      &= \bm{g}   &      &\text{ on }\Gamma_N, \\
  \bm{u}                       &= \bm{u}_D &      &\text{ on }\Gamma_D.
  \label{eq:elasticity-last-equation}
\end{alignat}
The constitutive equation for linear elastostatics
is given by the generalized Hooke's law
\begin{equation}
    \label{eq:hooke}
    \bm{\sigma}(\bm{u}) = 2 \mu \bm{\epsilon} (\bm{u}) + \lambda \mathrm{tr} \bm{\epsilon} (\bm{u}),
\end{equation}
where the linearized strain
\begin{equation}
    \bm{\epsilon} = \frac{1}{2} \left ( \nabla \bm{u} + \nabla \bm{u}^T \right ),
\end{equation}
is the symmetric part of the displacement gradient $\nabla \bm{u}$.
The first and second Lam\'e parameters, $\lambda$ and $\mu$,
are related to Young's modulus and Poisson's ratio, \(E\) and \(\nu\), by
\begin{equation}
  \lambda = \frac{E \nu}{(1 + \nu) (1 - 2\nu)}
  \quad\text{ and }\quad
  \mu = \frac{E}{2(1 + \nu)}.
\end{equation}

\subsection{DC PSE derivative operators}

Before describing their application to the proposed recovery method
we first summarize the derivation and computation of the discretization corrected (DC)
particle strength exchange (PSE) operators that we will use to compute derivatives \citep{reboux_etal_2012_selforganizing,
schrader_etal_2012_choosing}. The DC PSE method is similar to operators used to approximate derivatives, such as the finite difference method, corrected smoothed-particle hydrodynamics (SPH) method, reproducing kernel particle method (RKPM) and moving least squares (MLS) \citep{schrader_2011_discretizationcorrected}.
Their main advantage is that,
due to their meshless nature,
they can be applied relatively easily to arbitrarily complex geometries.
In this section,
we provide a succinct summary of the method
as background for the next section where we apply DC PSE
to compute the displacement gradient.
Further details on the derivation and numerical aspects of the DC PSE operators
can be found in the earlier contributions \citep{schrader_2011_discretizationcorrected,
reboux_etal_2012_selforganizing,
schrader_etal_2012_choosing}.

We define the DC PSE operator for approximating the spatial derivative $\boldsymbol{D}^{m,n}f(\bm{x})$ as:
\begin{equation}
    \label{dcpse_op}
	\boldsymbol{Q}^{m,n}f(\boldsymbol{x}_p)=\frac{1}{\epsilon(\boldsymbol{x}_p)^{m+n}}\sum_{\boldsymbol{x}_q \in N(\boldsymbol{x}_p)}\left(f(\boldsymbol{x}_q)\pm f(\boldsymbol{x}_p)\right)\frac{1}{\epsilon(\boldsymbol{x}_p)^2}\eta\left(\frac{\boldsymbol{x}_p-\boldsymbol{x}_q}{\epsilon(\boldsymbol{x}_p)}\right),
\end{equation}
where $\epsilon(\boldsymbol{x})$ is a spatially dependent scaling or resolution function, $\eta(\boldsymbol{x},\epsilon(\boldsymbol{x}))$ is a kernel function normalized by the factor $\epsilon(\boldsymbol{x}_p)^{-2}$, and  $N(\boldsymbol{x}_p)$ is the set of points in the support of the kernel function. 

We construct the DC PSE operators such that when we decrease the average spacing between nodes, $h(\boldsymbol{x}_p)\rightarrow 0$, the operator converges to the spatial derivative $\boldsymbol{D}^{m,n}f(\boldsymbol{x}_p)$ with an asymptotic rate $r$:
\begin{equation}
    \label{dcpse_cond}
	\boldsymbol{Q}^{m,n}f(\boldsymbol{x}_p) = \boldsymbol{D}^{m,n}f(\boldsymbol{x}_p) + \mathcal{O}(h(\boldsymbol{x}_p)^r),
\end{equation}
where  the component-wise average neighbor spacing is given as $h(\boldsymbol{x}_p) = \frac{1}{N}\sum\limits_{\boldsymbol{x}_q \in N(\boldsymbol{x})}|(x-x_q)|+|(y-y_q)|$, where $N$ is the number of nodes in the support of $\boldsymbol{x}_p$. Therefore, we need to find a kernel function $\eta(\boldsymbol{x})$ and a scaling relation $\epsilon(\boldsymbol{x}_p)$ as to satisfy Eq.~(\ref{dcpse_cond}). To achieve this, we replace the terms $f(\boldsymbol{x}_q)$ in Eq.~(\ref{dcpse_op}) with their Taylor series expansions around the point $\boldsymbol{x}_p$. This substitution gives:
\begin{equation}
\label{expansion}
\begin{split}
	\boldsymbol{Q}^{m,n}f(\boldsymbol{x}_p)&=\frac{1}{\epsilon(\boldsymbol{x}_p)^{m+n}}\sum_{\boldsymbol{x}_q \in N(\boldsymbol{x}_p)}\Biggl(\sum_{i=0}^{\infty}\sum_{j=0}^{\infty}\frac{(x_p-x_q)^i(y_p-y_q)^j(-1)^{i+j}}{i!j!}\boldsymbol{D}^{i,j}f(\boldsymbol{x}_p) \\ 
	&\pm  f(\boldsymbol{x}_p)\Biggr)\frac{1}{\epsilon(\boldsymbol{x}_p)^2}\eta\left(\frac{\boldsymbol{x}_p-\boldsymbol{x}_q}{\epsilon(\boldsymbol{x}_p)}\right). 
\end{split}
\end{equation}
It is convenient to re-write Eq.~(\ref{expansion}) in the form:
\begin{equation}
\label{expansion_moments}
\begin{split}
	\boldsymbol{Q}^{m,n}f(\boldsymbol{x}_p)&=\left(\sum_{i=0}^{\infty}\sum_{j=0}^{\infty}\frac{\epsilon(\boldsymbol{x}_p)^{i+j-m-n}(-1)^{i+j}}{i!j!}\boldsymbol{D}^{i,j}f(\boldsymbol{x}_p)Z^{i,j}(\boldsymbol{x}_p) \right) \\ 
	& \pm Z^{0,0}(\boldsymbol{x}_p) \epsilon(\boldsymbol{x}_p)^{-m-n}f(\boldsymbol{x}_p) ,
\end{split}
\end{equation}
where
\begin{equation}
	Z^{i,j}(\boldsymbol{x}_p) = \sum_{\boldsymbol{x}_q \in N(\boldsymbol{x}_p)}\frac{(x_p-x_q)^i(y_p-y_q)^j(-1)^{i+j}}{\epsilon(\boldsymbol{x}_p)^{i+j}}\frac{1}{\epsilon(\boldsymbol{x}_p)^2}\eta\left(\frac{\boldsymbol{x}_p-\boldsymbol{x}_q}{\epsilon(\boldsymbol{x}_p)}\right).
\end{equation}
We call $Z^{i,j}(\boldsymbol{x}_p)$ the {\it discrete moments}. Now, if we restrict the scaling parameter $\epsilon(\boldsymbol{x}_p)$ to converge at the same rate as the average spacing between points $h(\boldsymbol{x}_p)$, that is
\begin{equation}
	\frac{h(\boldsymbol{x}_p)}{\epsilon(\boldsymbol{x}_p)} = \text{constant}\label{scaling_rel},
\end{equation}
we find that the discrete moments $Z^{i,j}$ are $\mathcal{O}(1)$ as $h(\boldsymbol{x}_p)\rightarrow 0$.
Given Eq.~(\ref{scaling_rel}), the convergence behavior of Eq.~(\ref{dcpse_cond}) is determined by the coefficients of the terms $\epsilon(\boldsymbol{x}_p)^{i+j-m-n}\boldsymbol{D}^{i,j}$. We make the DC PSE operator satisfy Eq.~(\ref{dcpse_cond}) by ensuring the coefficients of the $(m,n)$ term is 1, and all other coefficients of $\epsilon(\boldsymbol{x}_p)^{i+j-m-n}\boldsymbol{D}^{i,j}$ for $i+j-m-n<r$ are 0. This results in the following set of conditions for the discrete moments,
\begin{equation}
	Z^{i,j}(\boldsymbol{x}_p) =
	\begin{cases}
		i!j!(-1)^{i+j} & i=m,j=n,\\
		0 & \alpha_{\text{min}} < i + j < r + m + n, \\
		< \infty & \text{otherwise},
	\end{cases} \label{moment_conditions}
\end{equation}
where $\alpha_{\text{min}}$ is 1 if $i+j$ is odd and 0 if even. This is due to the zeroth moment $Z^{0,0}$ canceling out for odd $i+j$. The choice of the factor $\epsilon(\boldsymbol{x}_p)^{-m-n}$ in Eq.~(\ref{dcpse_op}) simply acts as to simplify the expression of the moment conditions.

For the kernel function $\eta(\boldsymbol{x})$ to satisfy the $l$ conditions given in Eq.~(\ref{moment_conditions}) for arbitrary neighborhood node distributions, the operator must have $l$ degrees of freedom. This leads to the requirement that the support $N(\boldsymbol{x})$  of the kernel function has to include at least $l$ nodes. In this paper we use kernel functions for the operator of the form
\begin{equation}
	\eta(\boldsymbol{x}) = 
	\begin{cases}
		\sum\limits_{i,j}^{i+j<r+m+n}a_{i,j}x^iy^i=j e^{-x^2-y^2} & \sqrt{x^2+y^2}<r_c, \\
		0 & \text{otherwise.}
	\end{cases}\label{kernel_function}
\end{equation}
This is a monomial basis multiplied by an exponential window function, where $r_c$ sets the kernel support and the $a_{i,j}$ are scalars to be determined to satisfy the moment conditions in Eq.~(\ref{moment_conditions}). The cut-off radius $r_c$ should be set to include at least $l$ collocation nodes in the support $N(x)$. To construct the operator $\boldsymbol{Q}^{m,n}f(\boldsymbol{x}_p)$ at node $\boldsymbol{x}_p$, the coefficients are found by solving a linear system of equations from Eqs.~(\ref{kernel_function}) and (\ref{moment_conditions}). With our choice of kernel function we have,
\begin{equation}
	\boldsymbol{Q}^{m,n}f(\boldsymbol{x}_p) = \frac{1}{\epsilon(\boldsymbol{x}_p)^{m+n}}\sum_{\boldsymbol{x}_q \in N(\boldsymbol{x}_p)}\left(f(\boldsymbol{x}_q)\pm f(\boldsymbol{x}_p)\right)\boldsymbol{p}\left(\frac{\boldsymbol{x}_q-\boldsymbol{x}_p}{\epsilon(\boldsymbol{x}_p)}\right)\boldsymbol{a}^T(\boldsymbol{x}_p) e^{\frac{-(x_p-x_q)^2-(y_p-y_q)^2}{\epsilon(\boldsymbol{x}_p)}},  \label{dcpse_construct}
\end{equation}
where $\boldsymbol{p}(\boldsymbol{x})=\{p_1(x),p_2(x),\dots,p_l(x)\}$ and $\boldsymbol{a}(\boldsymbol{x})$ are vectors of the terms in the monomial basis and of their coefficients in Eq.~(\ref{kernel_function}), respectively. The normalization factor $\epsilon(\boldsymbol{x}_p)^{-2}$ is absent, since it cancels out through the coefficients $\boldsymbol{a}(\boldsymbol{x}_p)$. Using this formulation, the operator system becomes eassy to obtain.

The linear system for the kernel coefficients then is:
\begin{equation}
	\boldsymbol{A}(\boldsymbol{x}_p)\boldsymbol{a}^T(\boldsymbol{x}_p) = \boldsymbol{b},  \label{linear_system}
\end{equation}
where
\begin{align}
	\boldsymbol{A}(\boldsymbol{x}_p) &= 	\boldsymbol{B}(\boldsymbol{x}_p)^T\boldsymbol{B}(\boldsymbol{x}_p) \in \mathbb{R}^{l\times l}\\
	\boldsymbol{B}(\boldsymbol{x}_p)) &= 	\boldsymbol{E}(\boldsymbol{x}_p)^T\boldsymbol{V}(\boldsymbol{x}_p) \in \mathbb{R}^{k\times l}\\
	\boldsymbol{b}&=(-1)^{m+n}\boldsymbol{D}^{m,n}\boldsymbol{p}(\boldsymbol{x})|_{x=0} \in \mathbb{R}^{l\times 1}. 
\end{align}
The scalar $k>l$ is the number of nodes in the support of the operator, $l$ the number of moment conditions to be satisfied, and $\boldsymbol{V}(\boldsymbol{x}_p)$ the Vandermonde matrix constructed from the monomial basis $\boldsymbol{p}(\boldsymbol{x}_p)$.
The diagonal matrix $\boldsymbol{E}(\boldsymbol{x}_p)$ contains the square roots of the values of the exponential window function at the neighboring nodes in the operator support.
Further, we define $\{\boldsymbol{z}_i\}_{i=1}^{k}=\{\boldsymbol{x}_p-\boldsymbol{x}_q\}_{\boldsymbol{x}_q\in N(\boldsymbol{x}_p)}$, the set of vectors pointing from all neighboring nodes in the support of $\boldsymbol{x}_p$ to $\boldsymbol{x}_q$. So then explicitly
\begin{align}
	\boldsymbol{V}(\boldsymbol{x}_p) &= 
	\begin{pmatrix}
		p_1(\frac{\boldsymbol{z}_1}{\epsilon(\boldsymbol{x})}) & p_2(\frac{\boldsymbol{z}_1}{\epsilon(\boldsymbol{x})}) & \cdots & p_l(\frac{\boldsymbol{z}_1}{\epsilon(\boldsymbol{x})}) \\
		p_1(\frac{\boldsymbol{z}_2}{\epsilon(\boldsymbol{x})}) & p_2(\frac{\boldsymbol{z}_2}{\epsilon(\boldsymbol{x})}) & \cdots & p_l(\frac{\boldsymbol{z}_2}{\epsilon(\boldsymbol{x})})\\
		\vdots & \vdots & \ddots & \vdots \\
		p_1(\frac{\boldsymbol{z}_k}{\epsilon(\boldsymbol{x})}) & p_2(\frac{\boldsymbol{z}_k}{\epsilon(\boldsymbol{x})}) & \cdots & p_l(\frac{\boldsymbol{z}_k}{\epsilon(\boldsymbol{x})})
	\end{pmatrix}
	\in \mathbb{R}^{k\times l}\\
	\boldsymbol{E}(\boldsymbol{x}_p) &= \text{diag} \left( \left\{e^{\frac{-|\boldsymbol{z}_i|^2}{2\epsilon(\boldsymbol{x})}} \right\}_{i=1}^{k}\right)\in \mathbb{R}^{k\times k}.
\end{align}

Once the matrix $\bm{A}$ is constructed at each node $\bm{x}_p$, the linear systems can be solved for the coefficients $\bm{a}$ used in the DC PSE operators at each node as in Eq.~\eqref{dcpse_construct}.
The matrix $\bm{A}$ only depends on the number of moment conditions \(l\) and the local distribution of nodes in the support domain. The invertibility of $\bm{A}$ depends entirely on that of the Vandermonde matrix \(\bm{V}\), due to $\bm{E}$ being a diagonal matrix with non-zero entries.

\subsection{Recovery by DC PSE}

The problem of recovering the gradient
(and derived quantities such as stress and strain)
can be considered as a typical curve/surface fitting problem
where a function is sampled at a discrete number of points
and we seek to determine the functions derivative at those points.
The stress can be considered as a derived quantity,
by which we mean that it is not a primary solution variable
that must be computed to obtain a solution in displacements,
but can be computed during post-processing of the results
after the displacement is known.

After the displacement solution has been computed using the finite element method
(or other suitable numerical method, or from experimental measurements)
the displacement field is known at the nodal points of the discretized spatial domain.
However,
we are interested in the stress field,
which is related to the strain.
The strain is related to the gradient of the displacement.
Hence, we require an accurate method for computing derivatives at specified points,
given the value of the function at those points.
The DC PSE operators are ideally suited to this task.

Once the displacement gradient $\nabla \bm{u}$ has been computed (using the displacement field and the DC PSE spatial derivatives defined in
Eq.~\eqref{dcpse_construct}),
the strain and stress can be obtained by substitution of $\nabla \bm{u}$
into the constitutive equation \eqref{eq:hooke}.
Other derived quantities such as principal stresses and
stress invariants can be easily computed.
For example, the von Mises stress is defined as
\begin{equation}
    \sigma_\mathrm{Mises} = \sqrt{\frac{3}{2} \bm{s}:\bm{s}},
\end{equation}
where $s$ is the deviatoric stress tensor
\begin{equation}
    \bm{s}(\bm{u}) = \bm{\sigma}(\bm{u}) - \frac{1}{3} \mathrm{tr} \, \bm{\sigma}(\bm{u}) \bm{I}.
\end{equation}

We note briefly that this procedure
applies only to path-independent materials
because only the final displacement is used
to compute the stress.
For path-dependent materials,
such as those exhibiting plastic behavior,
a more elaborate procedure may be conceived
whereby the stress is computed at each time increment
using the current displacements and state-variables
stored at the nodes
(the state variables are typically stored only at the integration points
because that is where stress is usually computed).

The salient advantage of the method described in this section is that
all the derived quantities can be obtained from the known displacement field
(obtained in the current study either analytically
or using the finite element method)
and the constitutive equation.
This makes it straightforward to implement the proposed method
as a post-processing tool in existing software
irrespective of the methods used to solve the equations.

\subsection{Quantitative measures of accuracy}

DC PSE was originally developed as a collocation method,
and, unlike the finite element method,
the function values and gradients are defined at the nodes only.
Similarly to finite difference methods,
we do not use shape functions to compute function values
at arbitrary locations within the domain
so the integral error measures
typically used for assessing the solution accuracy of FEM
are not applicable.
Instead, we use discrete error measures,
namely the maximum absolute (\(\infty\)-norm) error and
the normalized root mean square error (NRMSE)
that are typically used for assessing the accuracy of collocation methods
to compare the predictions of the proposed method with analytical and finite element solutions.
The NRMSE between
the estimated nodal values \(\hat{\bm{y}}\) and
the reference nodal values \(\bm{y}\) is defined as
\begin{equation}
  \label{eq:nrmse}
  \mathrm{NRMSE}(\bm{y},\hat{\bm{y}}) = \frac{\sqrt{\mathrm{MSE}(\bm{y},\hat{\bm{y}})}}{\bm{y}_\mathrm{max} - \bm{y}_\mathrm{min}}, \quad
  \mathrm{MSE}(\bm{y},\hat{\bm{y}}) = \frac{1}{n} \sum_{i=1}^n (\bm{y}_i - \hat{\bm{y}}_i)^2,
\end{equation}
where \(n\) is the total number of nodes.
The maximum absolute (\(\infty\)-norm) error is defined as
\begin{equation}
  \label{eq:linf}
  ||\bm{y} - \hat{\bm{y}}||_\infty = \mathrm{max}_{1 \le i \le n} | \bm{y}_i - \hat{\bm{y}}_i |.
\end{equation}

\section{Results}
\label{sec:examples}

In this section,
we demonstrate the accuracy of our proposed meshfree recovery method
for some benchmark problems that have analytical solutions,
before applying the method to a practical example
of computing stress within the wall of an abdominal aortic aneurysm (AAA).
We compare the proposed meshfree DC PSE method
with the following commonly used finite elements:
triangles and tetrahedra with linear (P1) and quadratic (P2) Lagrange bases; and
quadrilaterals and hexahedra with linear (Q1) and quadratic (Q2) tensor-product Lagrange bases.

\FloatBarrier
\subsection{Franke's function}

The following function,
introduced by \citet{franke_1979_critical},
is commonly used to evaluate the performance of meshfree methods.
It consists of two Gaussian peaks and a sharper Gaussian dip
superimposed on a surface sloping towards the first quadrant,
and is defined as
\begin{multline}
  \label{eq:franke}
  u(x,y)
  = 0.75\,\mathrm{exp}\left[-\frac{(9x-2)^2      + (9y-2)^2}{4} \right]
  + 0.75\,\mathrm{exp}\left[-\frac{(9x+1)^2}{49} + \frac{(9y+1)^2}{10} \right] \\
  + 0.5 \,\mathrm{exp}\left[-\frac{(9x-7)^2      + (9y-3)^2)}{4} \right]
  - 0.2 \,\mathrm{exp}\left[-(9x-4)^2      + (9y-7)^2 \right].
\end{multline}

Given the solution \eqref{eq:franke},
we computed the derivatives numerically at the nodal points
using the FE and DC PSE methods
to compare the accuracy of the recovered gradients
\(\mathrm{\partial}u/\mathrm{\partial}x\) and \(\mathrm{\partial}u/\mathrm{\partial}y\)
against the analytical derivatives.
We used structured quadrilateral and triangular meshes,
as well as unstructured triangular meshes,
to discretize the spatial domain;
and used the same meshes for the FE and DC PSE methods
(disregarding the nodal connectivity when using DC PSE).
For the FE method, we used
linear (P1) and quadratic (P2) triangular, and
linear (Q1) and quadratic (Q2) quadrilateral finite elements.

\begin{figure}
  \centering
  \hfill
  \begin{subfigure}[b]{0.3\textwidth}
    \includegraphics[width=\textwidth]{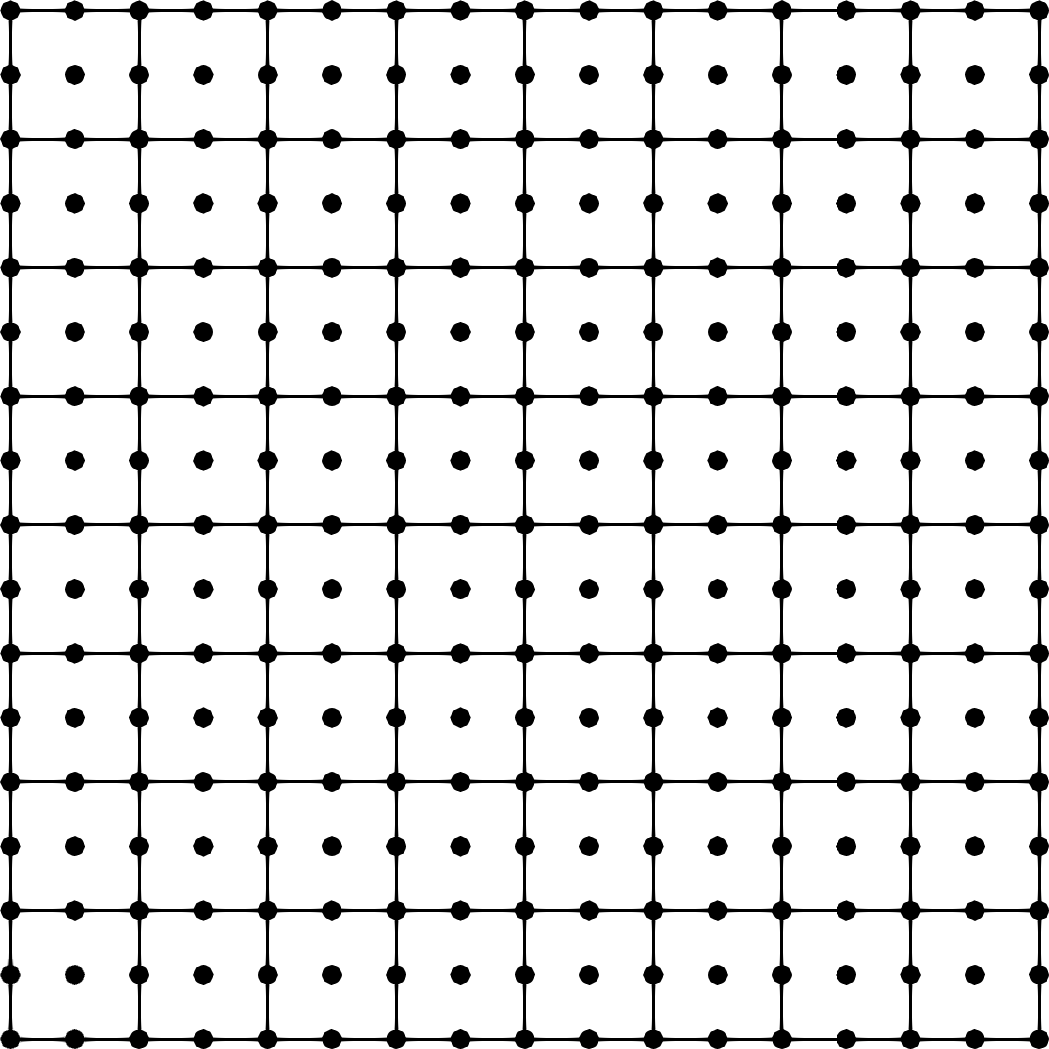}
    \caption{}
  \end{subfigure}
  \hfill
  \begin{subfigure}[b]{0.3\textwidth}
    \includegraphics[width=\textwidth]{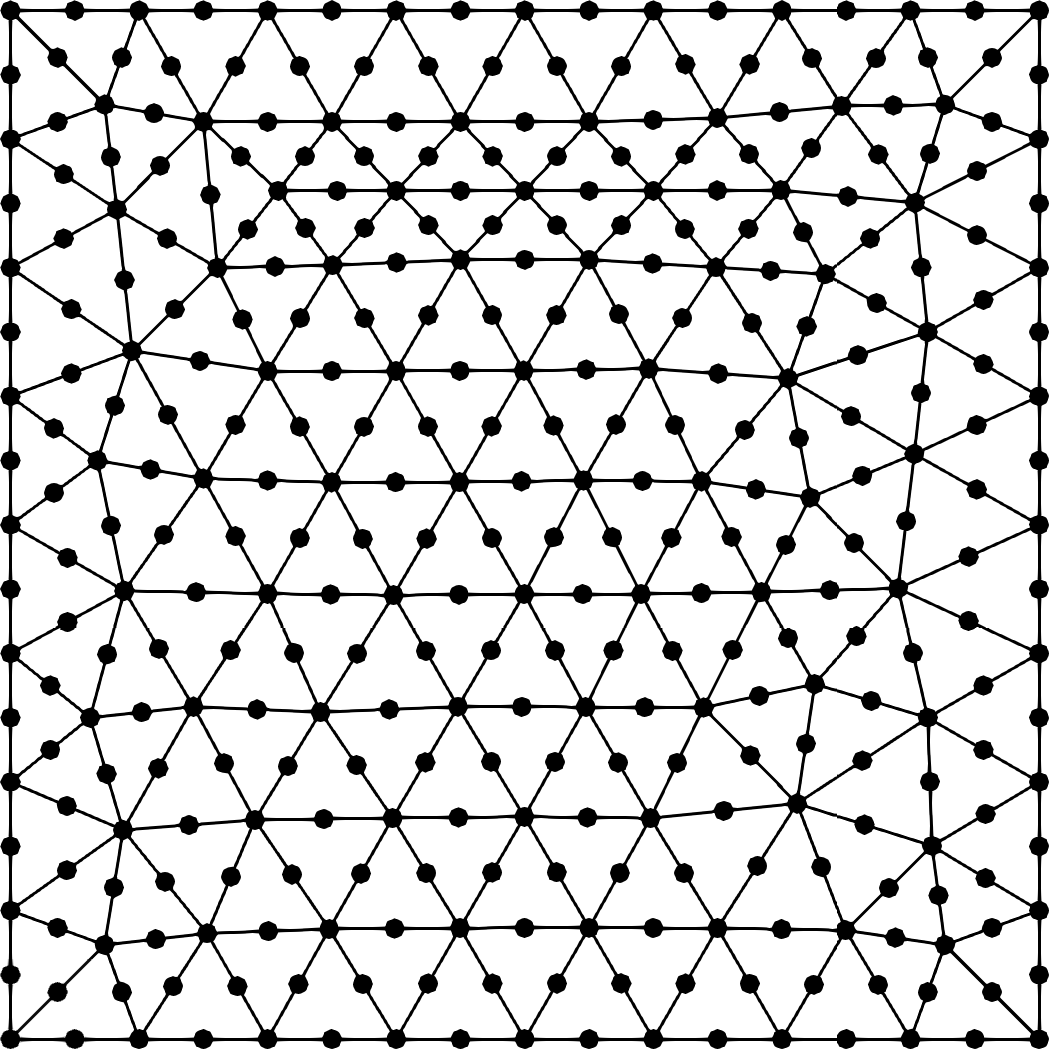}
    \caption{}
  \end{subfigure}
  \hfill{}
  \caption{
    The
    (a) coarsest structured quadratic quadrilateral (Q2) and
    (b) coarsest unstructured quadratic triangle (P2)
    meshes used for Franke's function benchmark problem.
    The meshes have
    (a) 289 nodes (\(h=0.063\)) and 64 elements, and
    (b) 357 nodes (\(h=0.056\)) and 162 elements, respectively.
    The coarsest linear quadrilateral and triangle meshes
    have a similar number of nodes with more elements.
  }
  \label{fig:franke-meshes}
\end{figure}

Fig.~\ref{fig:franke-meshes}
shows representative quadrilateral and triangle meshes used for this problem.
The nodal distribution corresponding to the quadrilateral elements
for a given Q2 mesh and a refined Q1 mesh
obtained by splitting the elements of the Q1 mesh are equivalent,
with both meshes having the same number of nodes, node spacing \(h\)
and nodal coordinates
(only the number of elements and element nodal connectivity is different).
Hence, the DC PSE derivatives for these two grids will be equivalent
when the function values at the nodes are the same.
The triangular P1 and P2 meshes share similar similarities.

\begin{figure}
  \centering
  \begin{subfigure}[b]{0.45\textwidth}
    \includegraphics[width=\textwidth]{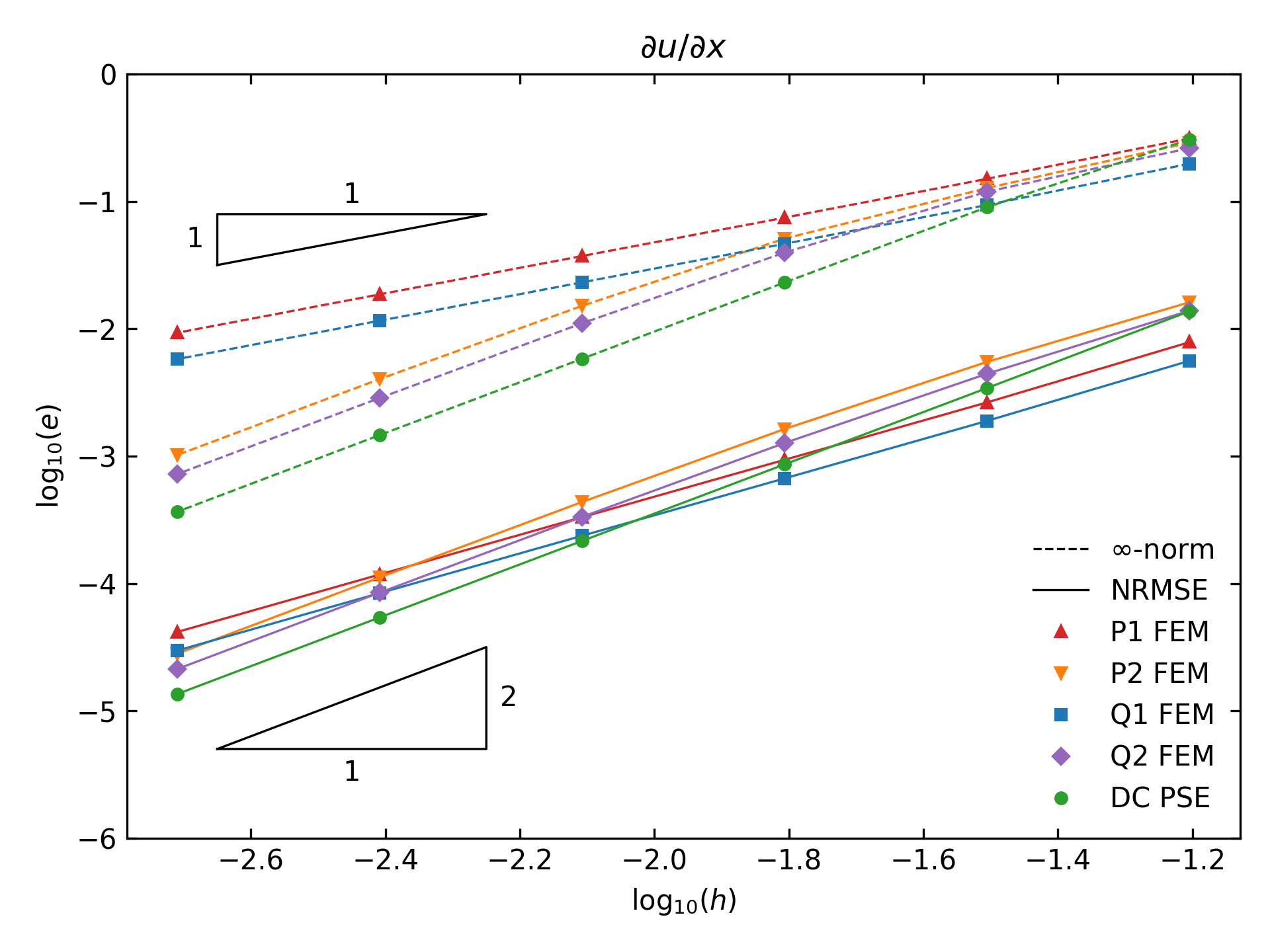}
    \caption{}
  \end{subfigure}
  \hfill
  \begin{subfigure}[b]{0.45\textwidth}
    \includegraphics[width=\textwidth]{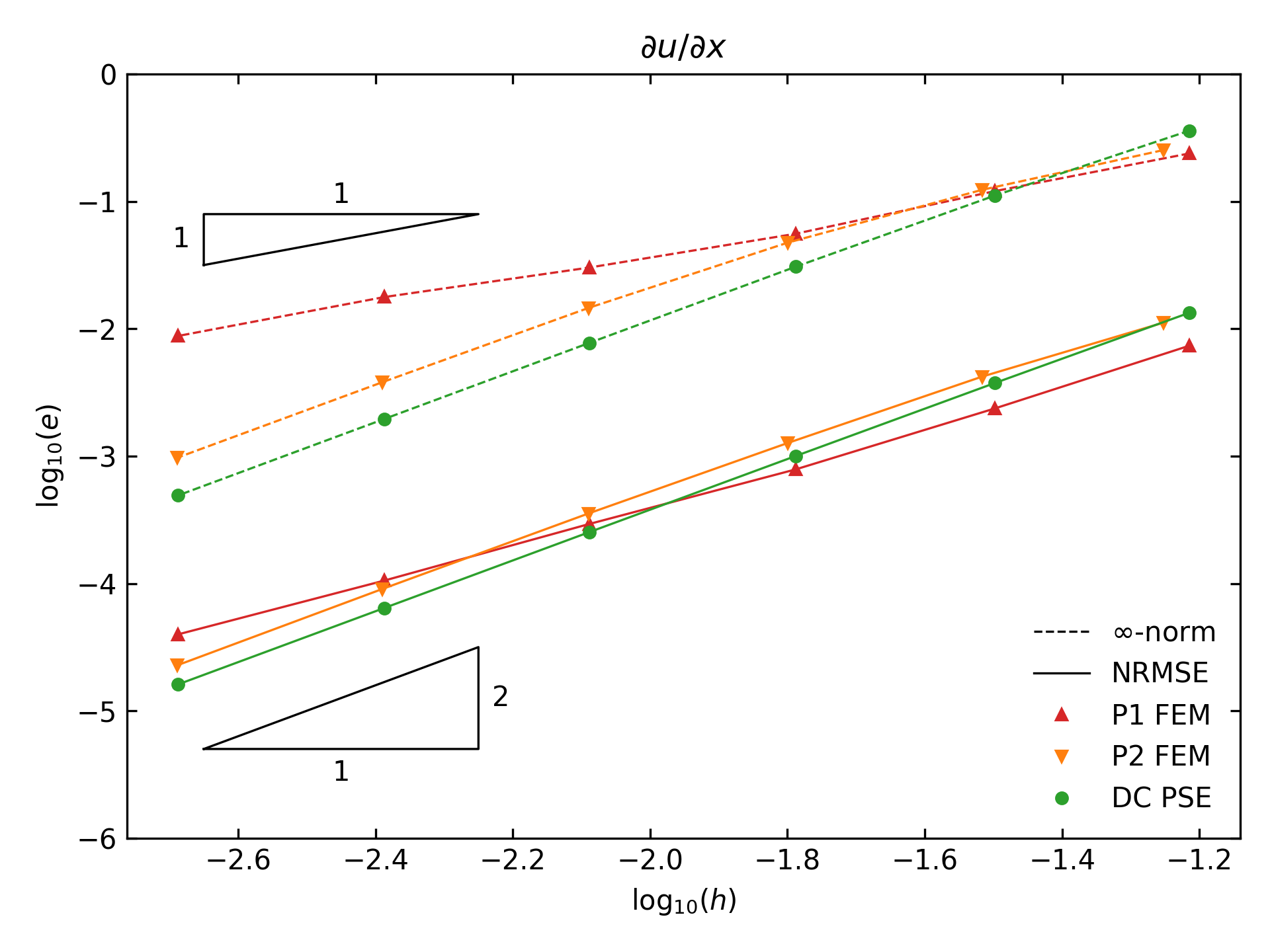}
    \caption{}
  \end{subfigure}
  \\
  \begin{subfigure}[b]{0.45\textwidth}
    \includegraphics[width=\textwidth]{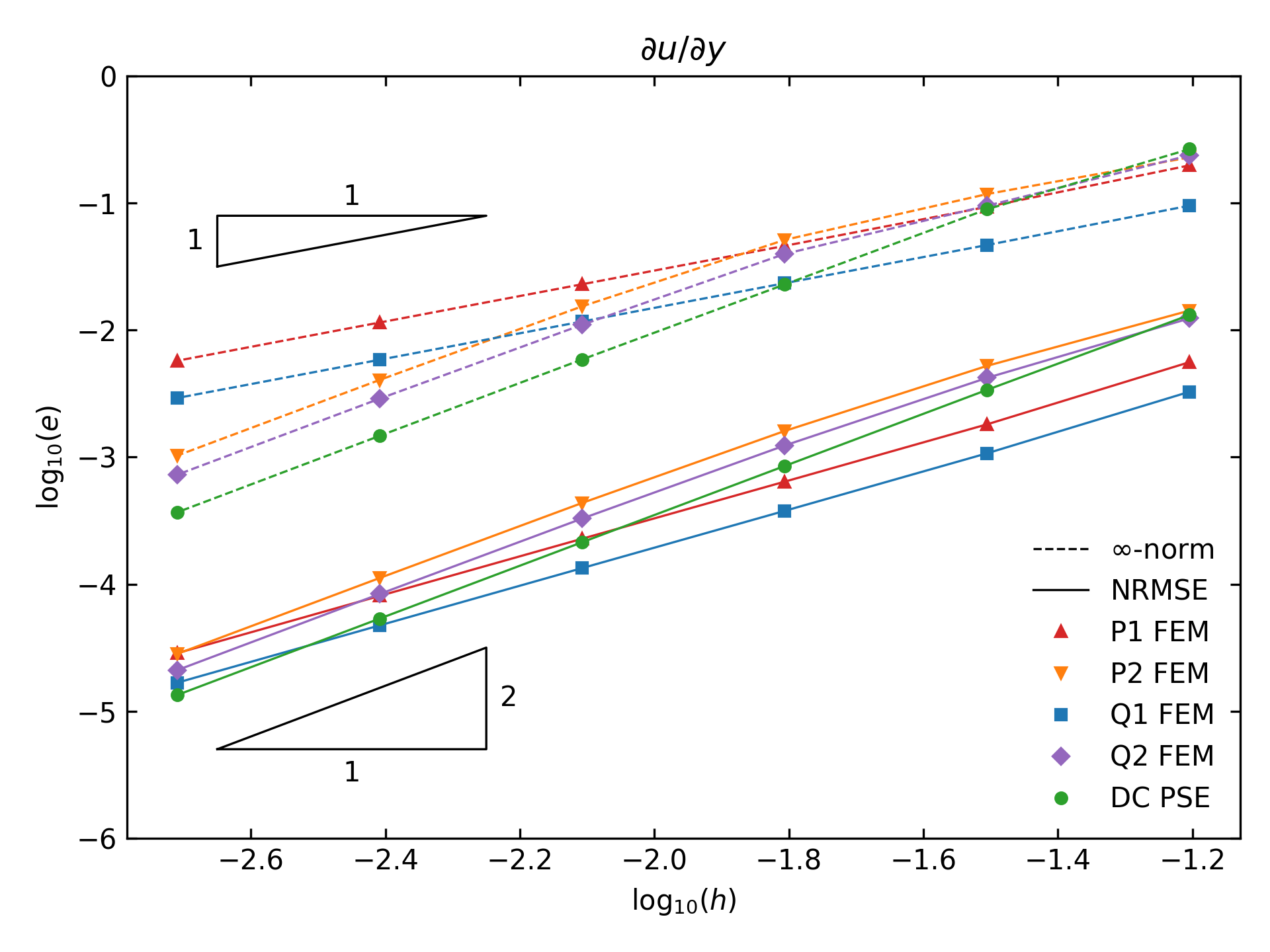}
    \caption{}
  \end{subfigure}
  \hfill
  \begin{subfigure}[b]{0.45\textwidth}
    \includegraphics[width=\textwidth]{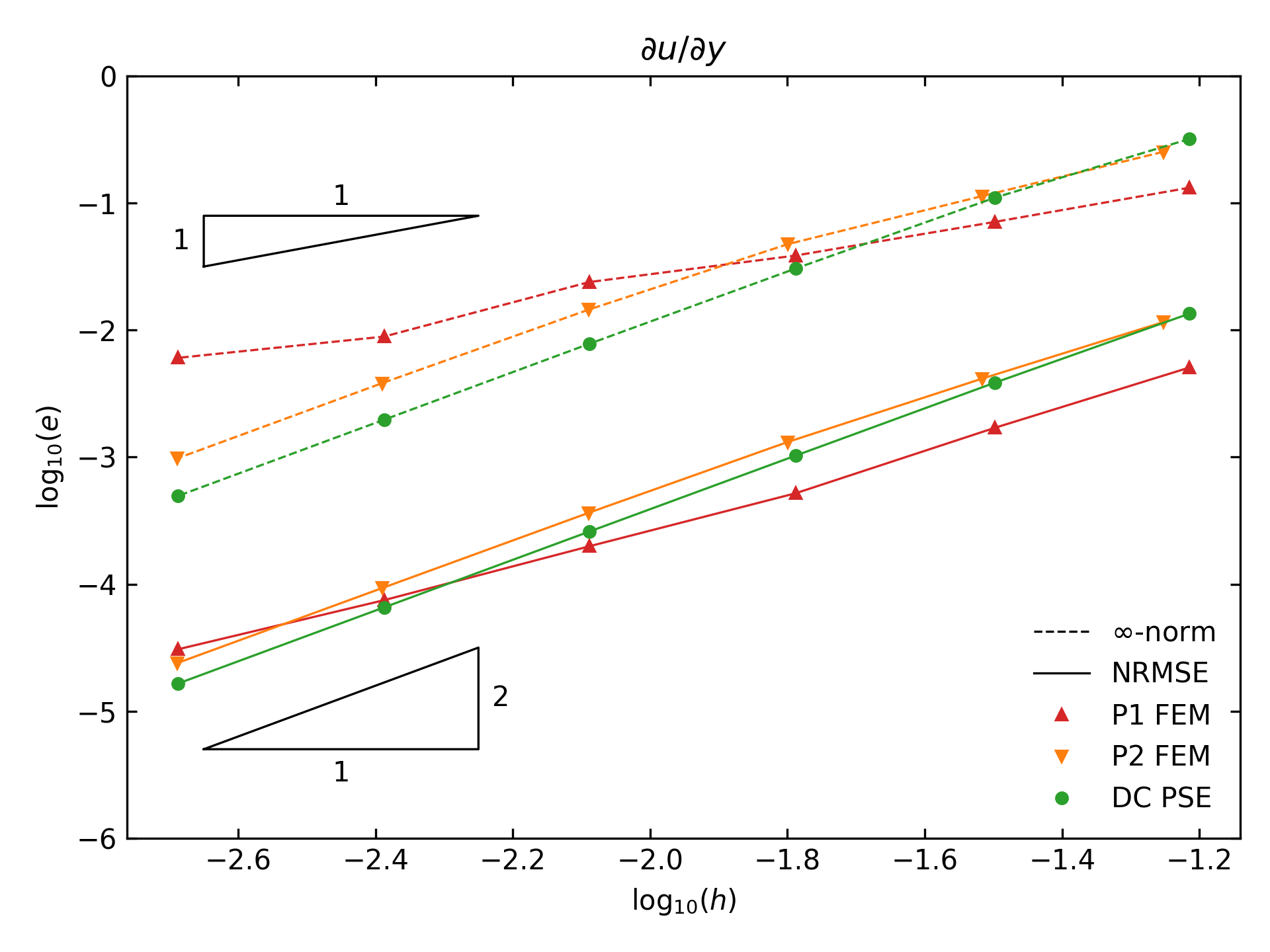}
    \caption{}
  \end{subfigure}
  \caption{
    Convergence of numerical derivatives
    (a, b) \(\mathrm{\partial}u/\mathrm{\partial}x\) and 
    (c, d) \(\mathrm{\partial}u/\mathrm{\partial}y\)
    of Franke's function
    using
    (a, c) structured triangle
    and structured quadrilateral meshes and
    (b, d) unstructured triangle meshes.
    The derivatives were computed using 
    linear (P1) and quadratic (P2) triangle finite elements,
    linear (Q1) and quadratic (Q2) quadrialteral finite elements,
    and the proposed meshfree DC PSE based method.
    The normalized node spacing \(h\) for a 2D mesh with \(n\) nodes
    is defined as \(h = 1/(n^{1/2}-1)\).
  }
  \label{fig:franke-convergence}
\end{figure}

Fig.~\ref{fig:franke-convergence}
shows the convergence of the normalized root mean square error (NRMSE)
and the maximum absolute ($\infty$-norm) error
in the derivatives
\(\mathrm{\partial}u/\mathrm{\partial}x\) and \(\mathrm{\partial}u/\mathrm{\partial}y\)
of Franke's function
computed using the FE and DC PSE methods.
As expected, the derivatives computed using quadratic (P2 and Q2) finite elements and DC PSE
converge more rapidly,
and are significantly more accurate,
than those computed using linear (P1 and Q1) finite elements.
Moreover, the DC PSE derivatives are slightly more accurate
and converge slightly faster than the quadratic finite elements
for both the quadrilateral and triangle meshes.
These results suggest that DC PSE is well suited
for accurately recovering the gradients of a given function.

\FloatBarrier
\subsection{Infinite plate with a circular hole}

The axially loaded infinite plate with a circular hole
is a commonly used benchmark problem for stress analysis
\citep{
  kelly_etal_1983_posteriori,
  zienkiewicz_zhu_1992_superconvergent_part1_recovery,
  boroomand_zienkiewicz_1997_recovery}.
The elastic stresses around the circumference of the hole are of primary interest
and can be obtained analytically.
The problem is applicable to infinitely thin structures such as a steel plate
(plane stress)
and infinitely thick structures such as a borehole
(plane strain).
The analytical solution for the planar stress
is valid for both plane strain and plane stress.
Here, we consider only plane strain loading conditions,
where the displacement perpendicular to the plane of loading is assumed to be zero,
because the stress recovery method requires the full displacement field
for computing the strain and stress numerically.

The material is linear elastic with
Young's modulus $E=200$ GPa and Poisson's ratio $\nu=0.3$.
The analytical solution in polar coordinates \(r,\theta\)
for the displacement components---%
which applies only for plane strain---%
is given by
\citep{timoshenko_goodier_1951_theory} %
\begin{align}
  u_x(r,\theta) &= \frac{\sigma_0}{8\mu} \left(
                  \frac{r}{a}(k+1)\mathrm{cos}\theta
                  + \frac{2a}{r} ( (1+k)\mathrm{cos}\theta + \mathrm{cos}3\theta )
                  - \frac{2a^3}{r^3}\mathrm{cos}3\theta \right), \\
  u_y(r,\theta) &= \frac{\sigma_0}{8\mu} \left(
                  \frac{r}{a}(k-3)\mathrm{sin}\theta
                  + \frac{2a}{r} ( (1-k)\mathrm{sin}\theta + \mathrm{sin}3\theta )
                  - \frac{2a^3}{r^3}\mathrm{sin}3\theta \right),
\end{align}
where \(\sigma_0\) is the far field traction
and the Kolosov constant \(k=3-4\nu\) for plane strain,
with \(\mu=E/(2(1+\nu))\) and \(\nu\)
being Lam\'{e}'s second parameter (shear modulus)
and Poisson's ratio, respectively.
The analytical solution for the stress in the plane---%
valid for both plane strain and plane stress---%
is given by the Kirsch equations
\citep{timoshenko_goodier_1951_theory,kirsch_1898_theorie}%
\begin{align}
    \sigma_{xx}(r,\theta) &= \sigma_0 \left [ 1 - \frac{a^2}{r^2} \left ( \frac{3}{2} \mathrm{cos} \, 2 \theta + \mathrm{cos} \, 4 \theta \right ) + \frac{3a^4}{2r^4} \mathrm{cos} \, 4 \theta \right], \\
    \sigma_{yy}(r,\theta) &= \sigma_0 \left [   - \frac{a^2}{r^2} \left ( \frac{1}{2} \mathrm{cos} \, 2 \theta - \mathrm{cos} \, 4 \theta \right ) - \frac{3a^4}{2r^4} \mathrm{cos} \, 4 \theta \right], \\
    \sigma_{xy}(r,\theta) &= \sigma_0 \left [   - \frac{a^2}{r^2} \left ( \frac{1}{2} \mathrm{sin} \, 2 \theta + \mathrm{sin} \, 4 \theta \right ) + \frac{3a^4}{2r^4} \mathrm{sin} \, 4 \theta \right].
\end{align}

\begin{figure}
  \centering
  \hfill
  \begin{subfigure}[b]{0.3\textwidth}
    \includegraphics[width=\textwidth]{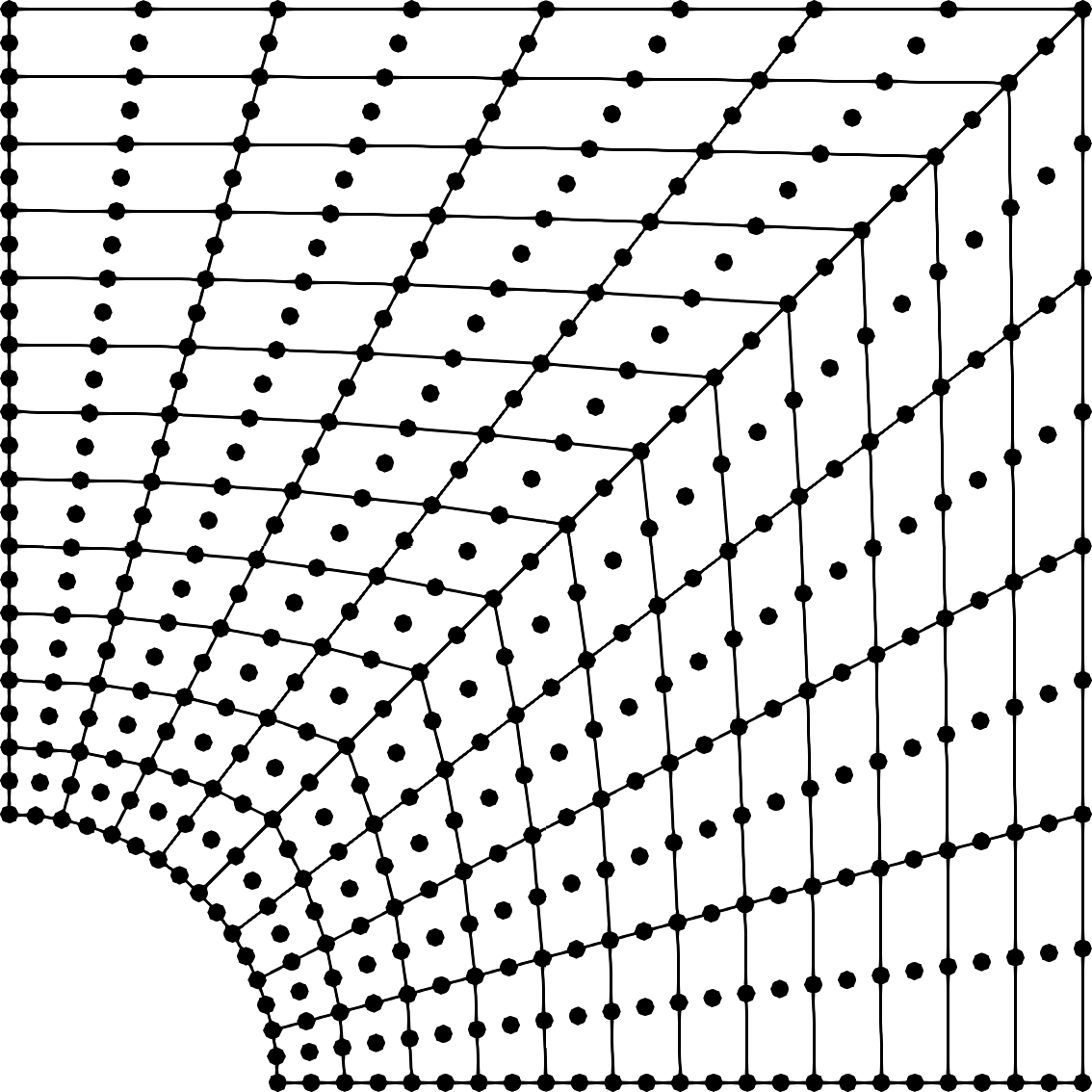}
    \caption{}
  \end{subfigure}
  \hfill
  \begin{subfigure}[b]{0.3\textwidth}
    \includegraphics[width=\textwidth]{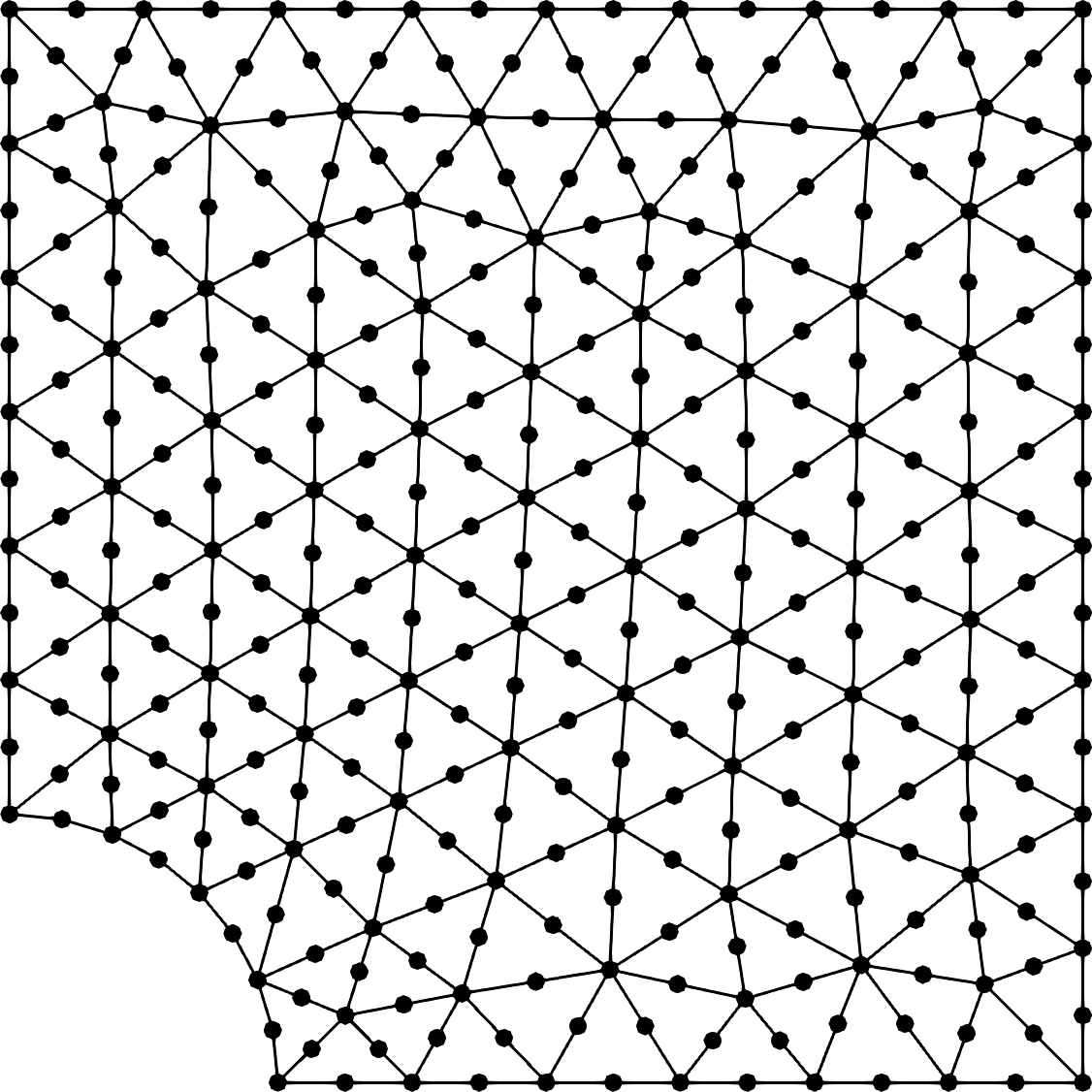}
    \caption{}
  \end{subfigure}
  \hfill{}
  \caption{
    The
    (a) coarsest structured quadratic quadrilateral (Q2) mesh
    with 425 nodes (\(h=0.051\)) and  96 elements, and
    (b) coarsest unstructured quadratic triangle (P2) mesh
    with 357 nodes (\(h=0.056\)) and 162 elements
    used for the infinite plate with hole benchmark problem.
    The coarsest linear quadrilateral and triangle meshes
    have a similar number of nodes with more elements.
  }
  \label{fig:platewithhole-meshes}
\end{figure}

\begin{figure}
  \centering
  \begin{subfigure}[b]{0.45\textwidth}
    \includegraphics[width=\textwidth]{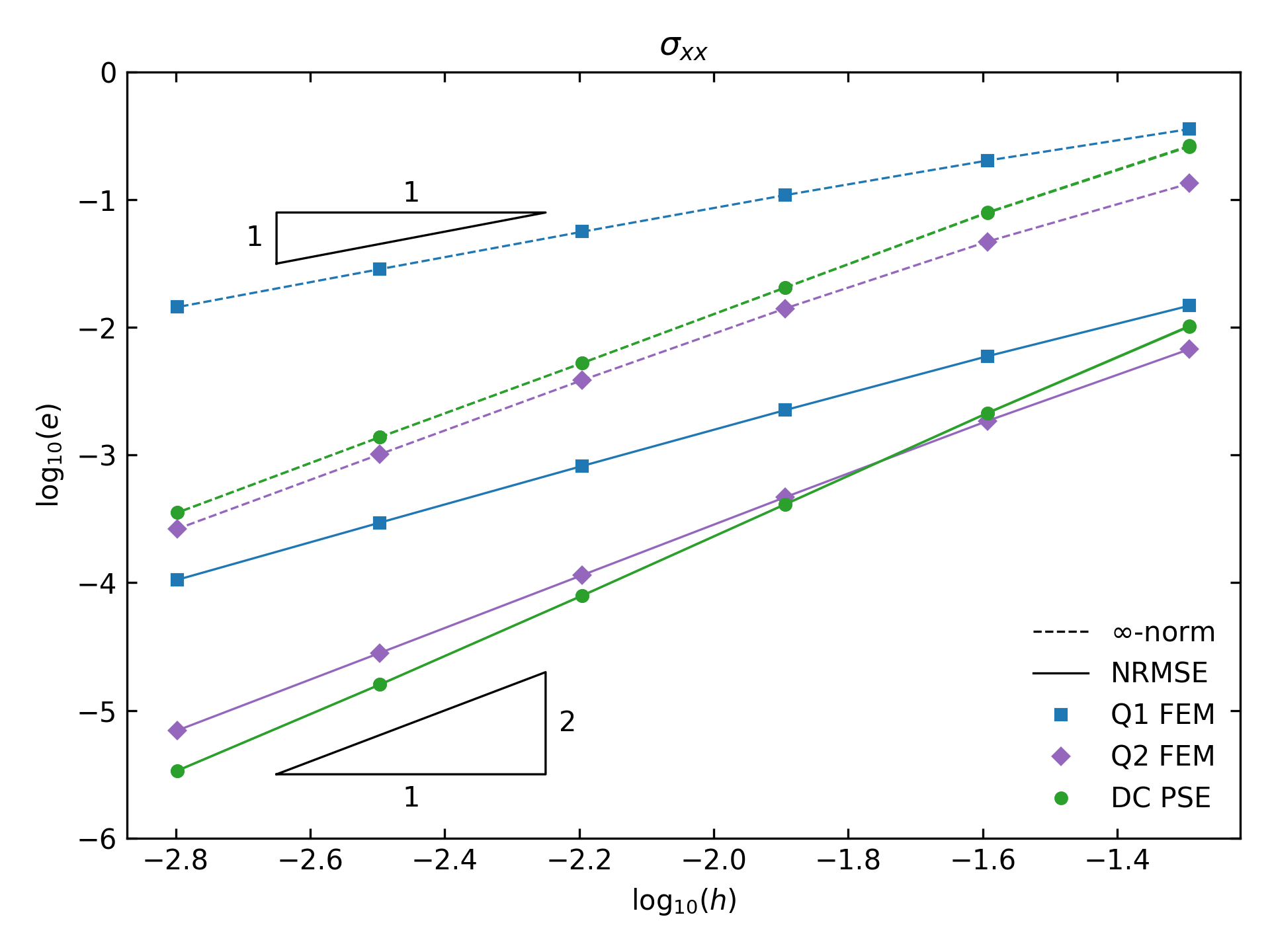}
    \caption{}
  \end{subfigure}
  \begin{subfigure}[b]{0.45\textwidth}
    \includegraphics[width=\textwidth]{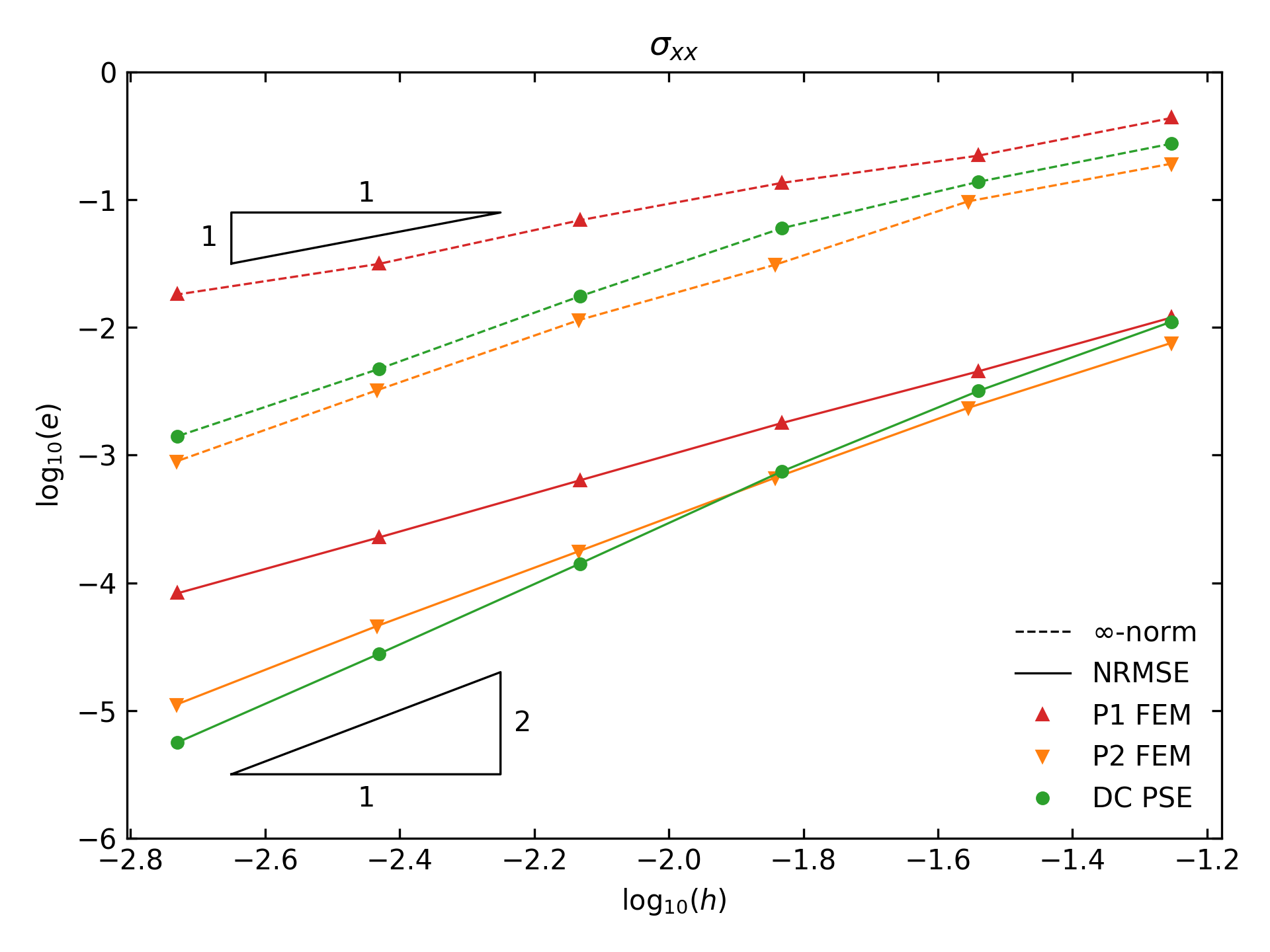}
    \caption{}
  \end{subfigure}  
  \\
  \begin{subfigure}[b]{0.45\textwidth}
    \includegraphics[width=\textwidth]{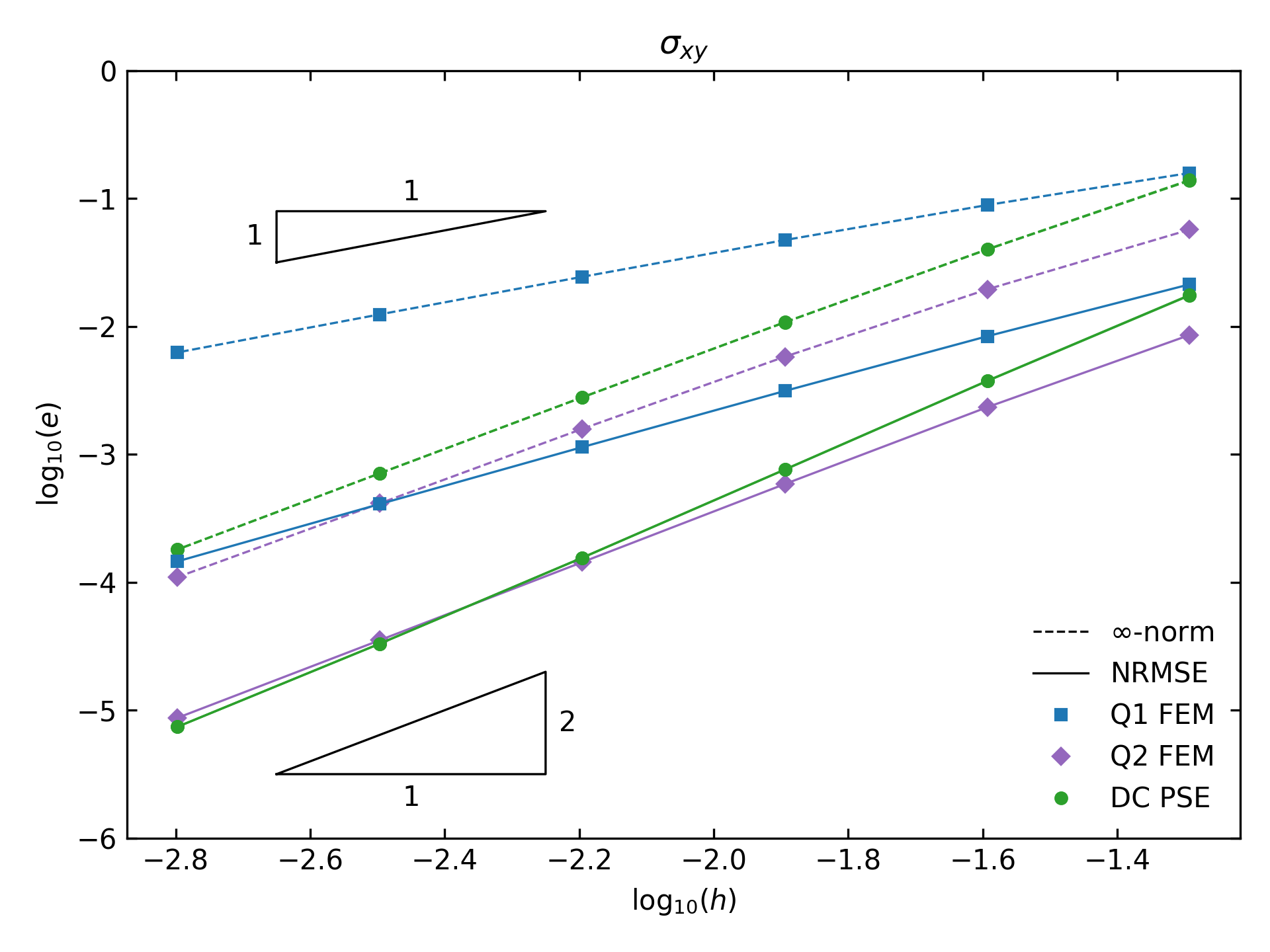}
    \caption{}
  \end{subfigure}
  \begin{subfigure}[b]{0.45\textwidth}
    \includegraphics[width=\textwidth]{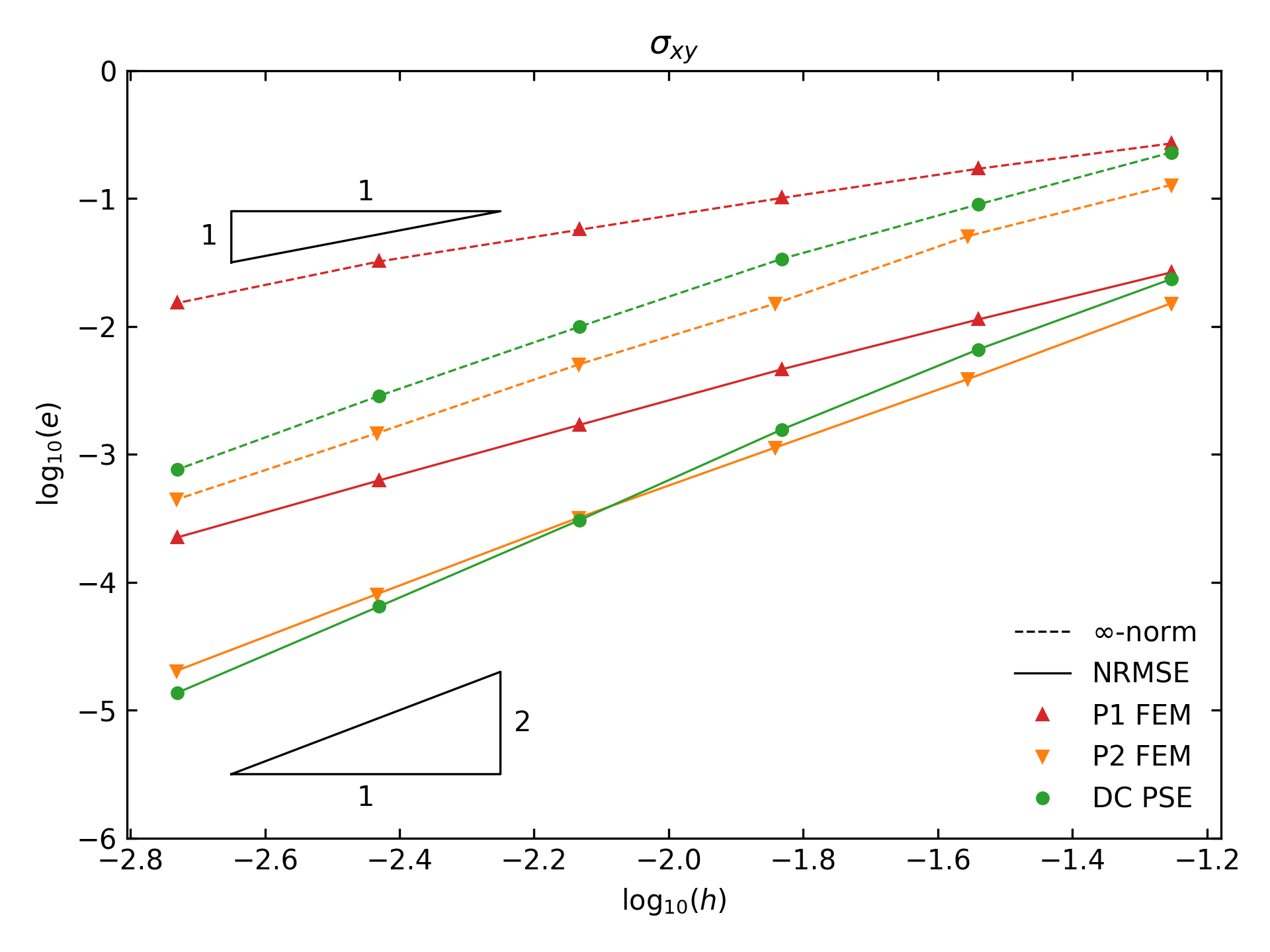}
    \caption{}
  \end{subfigure}  
  \\
  \begin{subfigure}[b]{0.45\textwidth}
    \includegraphics[width=\textwidth]{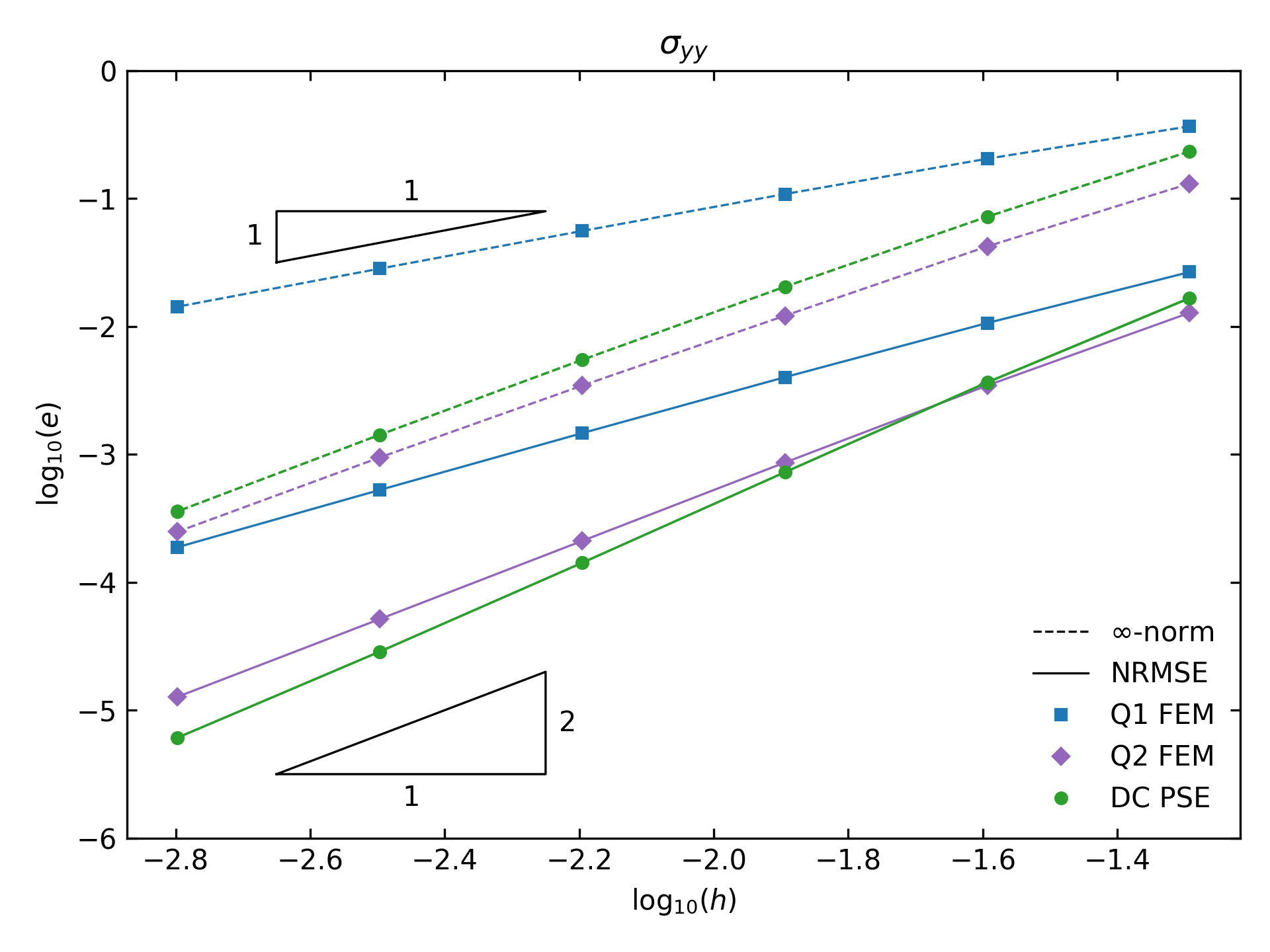}
    \caption{}
  \end{subfigure}
  \begin{subfigure}[b]{0.45\textwidth}
    \includegraphics[width=\textwidth]{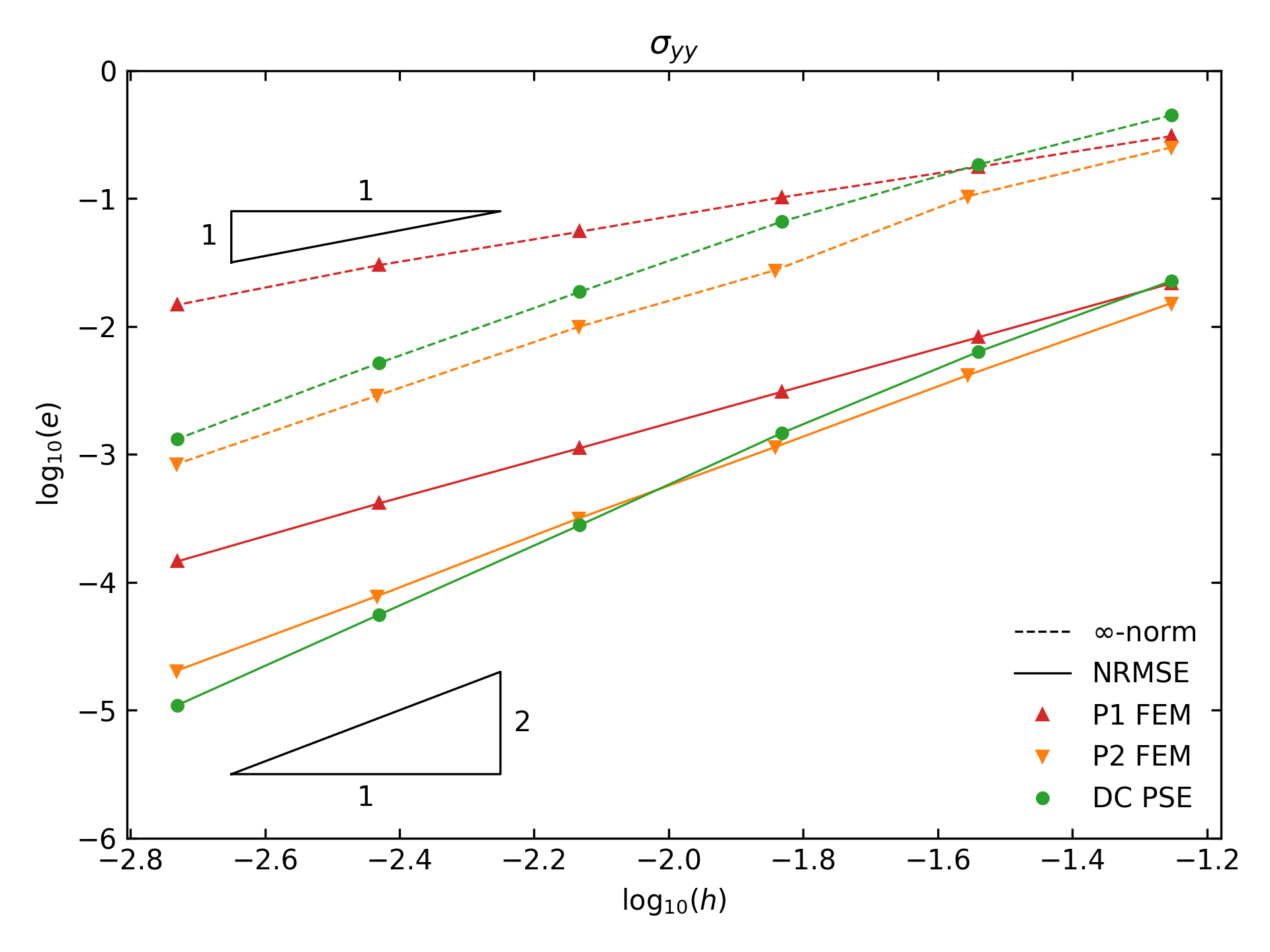}
    \caption{}
  \end{subfigure}  
  \caption{
    Convergence of stress components
    \(\sigma_{xx}\), \(\sigma_{xy}\) and \(\sigma_{yy}\)
    for the infinite plate with hole benchmark problem
    using
    (a, c, e) structured quadrilateral and
    (b, d, f) unstructured triangle meshes.
    The normalized node spacing \(h\) for a 2D mesh with \(n\) nodes
    is defined as \(h = 1/(n^{1/2}-1)\).
  }
  \label{fig:platewithhole-convergence-stress}
\end{figure}

Fig.~\ref{fig:platewithhole-meshes}
shows some representative meshes used for this problem.
Using the analytical solution for the displacements at the nodal points,
we computed the derivatives and derived quantities (strain and stress) numerically
using the FE method with nodal averaging
and the proposed DC PSE based meshfree recovery method
to compare the accuracy of the stress components.
Fig.~\ref{fig:platewithhole-convergence-stress}
shows the convergence of each stress component
computed at the nodal points.
As expected, we observe first and second order convergence
for the linear (P1 and Q1) and quadratic (P2 and Q2) finite elements, respectively.
The stress computed using DC PSE 
converges at the same rate,
and has similar accuracy,
as the quadratic finite elements.

\FloatBarrier
\subsection{3D cantilever beam}

As a benchmark 3D example,
we consider a rectangular beam with 
length $L = 10$ cm, width $2a = 2$ cm and height $2b = 2$ cm.
The beam is fixed (weakly) on the end $z=L$
and subjected to a transverse shear force $F = 1$ N
in the negative $y$-direction at the opposite end where $z=0$
\citep{bishop_2014_displacementbased}.
The material is linear elastic
with Young’s modulus $E = 10^7$ N/cm$^2$
and Poisson’s ratio $\nu = 0.3$. %
The analytical solution \citep{barber_2010_elasticity,bishop_2014_displacementbased}
is given as follows for the displacement field:
\begin{gather}
    u_x = -\frac{F \nu}{E I} x y z, \\
    u_y = \frac{F}{E I} \left[ \frac{\nu}{2} (x^2 - y^2)z - \frac{1}{6} z^3 \right], \\
    \begin{split}
    u_z = \frac{F}{E I} \left[ \frac{1}{2} y (\nu x^2 + z^2)
    + \frac{1}{6} \nu y^3 
    + (1 + \nu) \left( b^2 y - \frac{1}{3} y^3 \right)
    - \frac{1}{3} a^2 \nu y \right. \\
    \left. - \frac{4 a^3 \nu}{\pi^3} \sum_{n=1}^\infty \frac{(-1)^n}{n^3} \mathrm{cos}(n \pi x / a)
    \frac{\mathrm{sinh}(n \pi y / a)}{\mathrm{cosh}(n \pi b / a)}
    \right],
    \end{split}
\end{gather}
and for the stress field:
\begin{gather}
    \sigma_{xx} = \sigma_{xy} = \sigma_{yy} = 0, \\
    \sigma_{zz} = \frac{F}{I} y z, \\
    \sigma_{xz} = \frac{F}{I} \frac{2 a^2}{\pi^2} \frac{\nu}{1 + \nu}
      \sum_{n=1}^\infty \frac{(-1)^n}{n^2} \mathrm{sin}(n \pi x / a)
      \frac{\mathrm{sinh}(n \pi y / a)}{\mathrm{cosh}(n \pi b / a)}, \\
    \sigma_{yz} = \frac{F}{I} \frac{b^2 - y^2}{2} 
      + \frac{F}{I} \frac{\nu}{1 + \nu}
      \left[
      \frac{3x^2 - a^2}{6} - \frac{2a^2}{\pi^2}
      \sum_{n=1}^\infty \frac{(-1)^n}{n^2} \mathrm{cos}(n \pi x / a)
      \frac{\mathrm{cosh}(n \pi y / a)}{\mathrm{cosh}(n \pi b / a)}
      \right],
\end{gather}
where $I=4ab^3/3$ is the second moment of area about the $x$-axis.

\begin{figure}
  \centering
  \hfill
  \begin{subfigure}[b]{0.4\textwidth}
    \includegraphics[width=\textwidth]{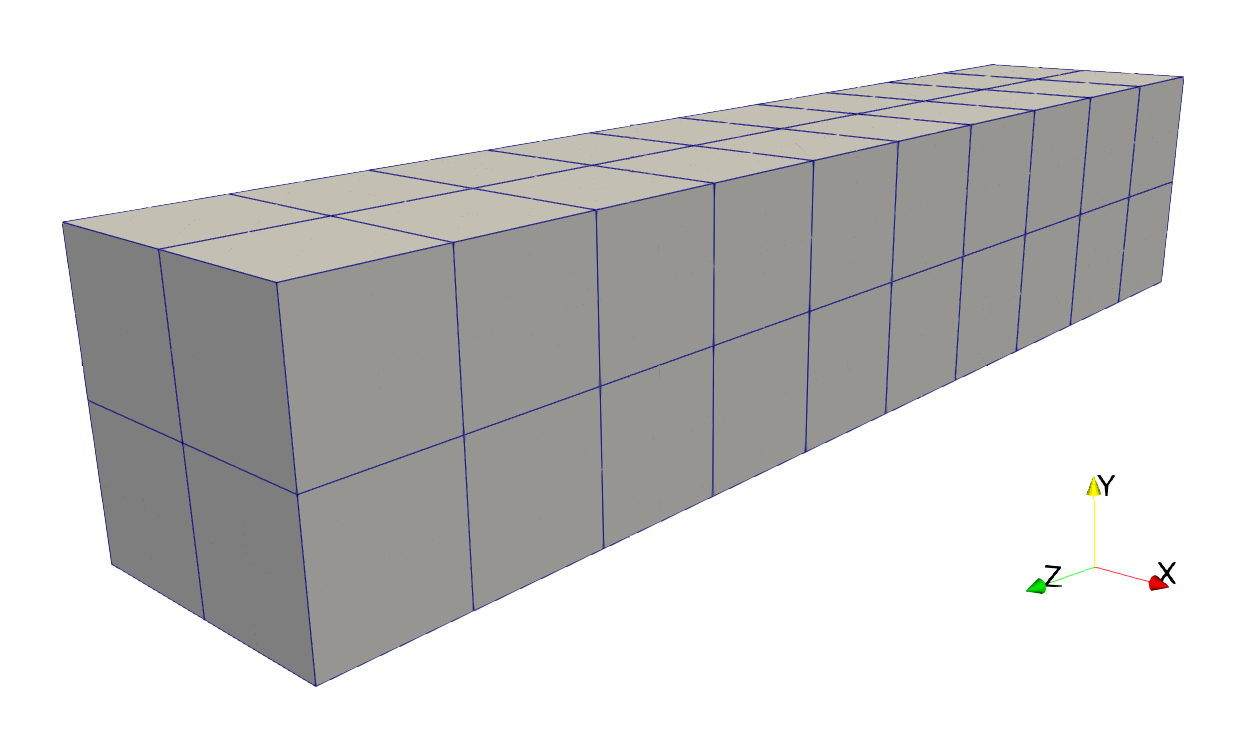}
    \caption{}
  \end{subfigure}
  \hfill
  \begin{subfigure}[b]{0.4\textwidth}
    \includegraphics[width=\textwidth]{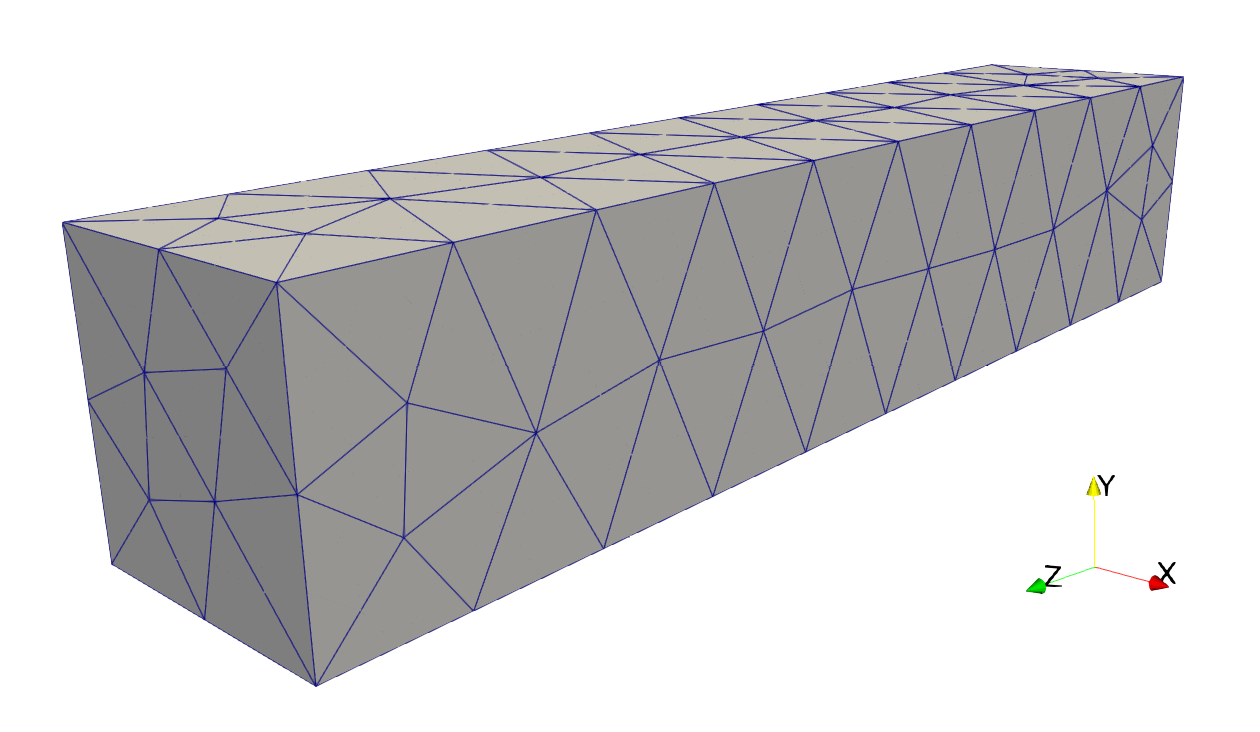}
    \caption{}
  \end{subfigure}
  \hfill{}
  \caption{
    The
    (a) coarsest structured quadratic hexahedral (Q2) mesh
    with 525 nodes (\(h=0.1415\)) and 40 elements, and
    (b) coarsest unstructured quadratic tetrahedral (P2) mesh
    with 587 nodes (\(h=0.1356\)) and 260 elements
    used for the 3D cantilever benchmark problem.
    The coarsest linear hexahedral and tetrahedral meshes
    have a similar number of nodes with more elements.
  }
  \label{fig:cantilever-meshes}
\end{figure}

\begin{figure}
  \centering
  \begin{subfigure}[b]{0.45\textwidth}
    \includegraphics[width=\textwidth]{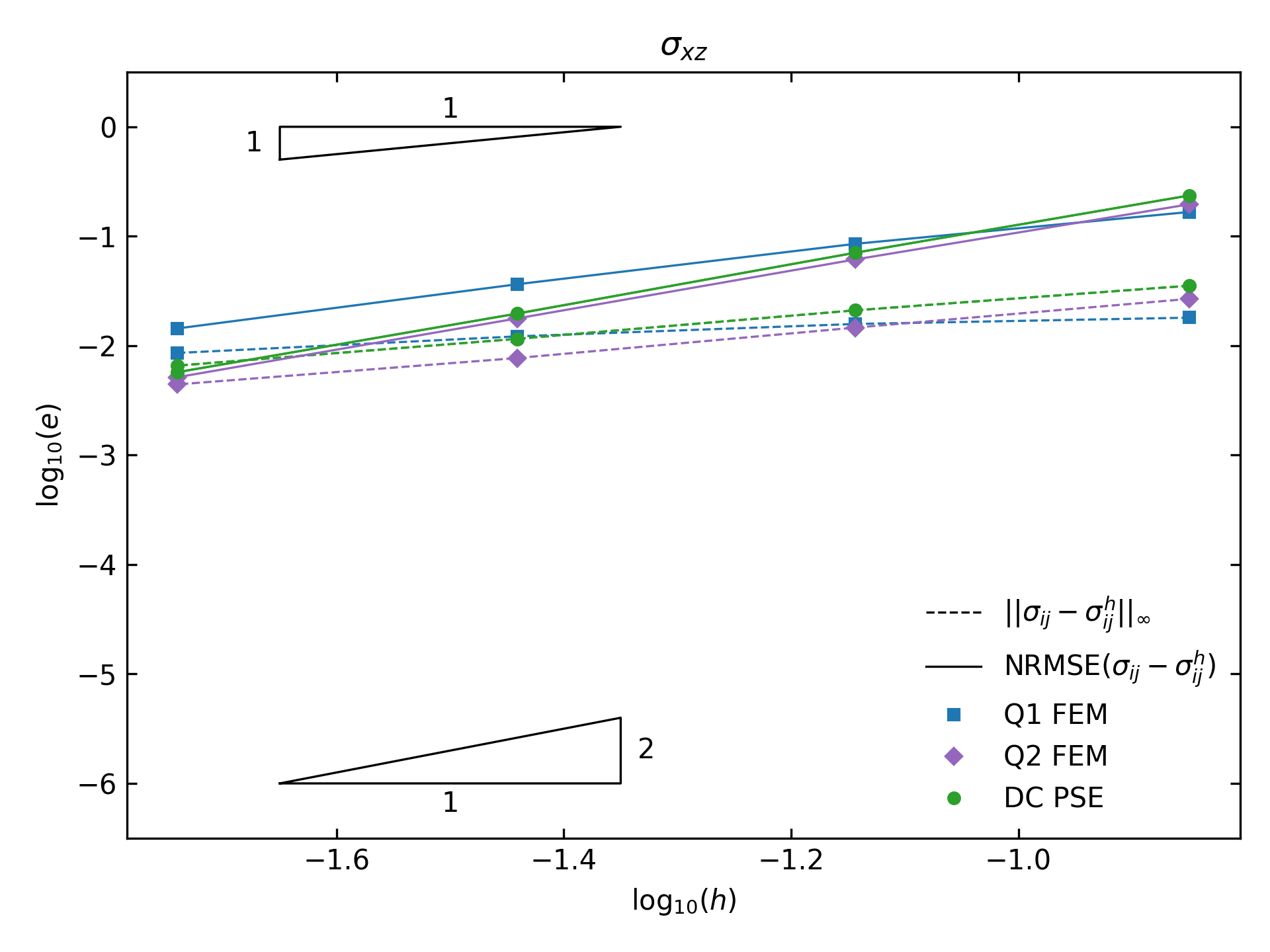}
    \caption{}
  \end{subfigure}
  \begin{subfigure}[b]{0.45\textwidth}
    \includegraphics[width=\textwidth]{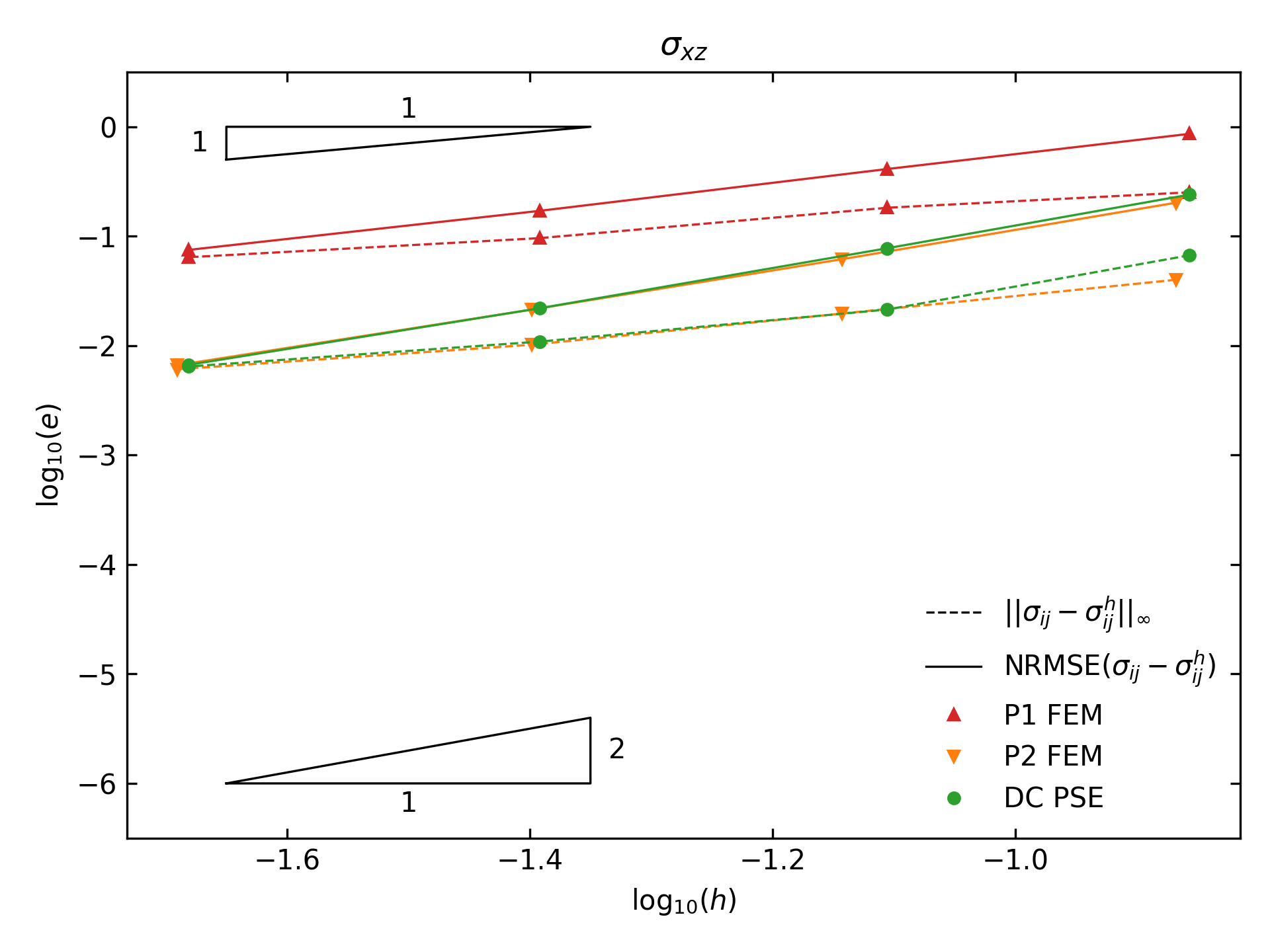}
    \caption{}
  \end{subfigure}  
  \\
  \begin{subfigure}[b]{0.45\textwidth}
    \includegraphics[width=\textwidth]{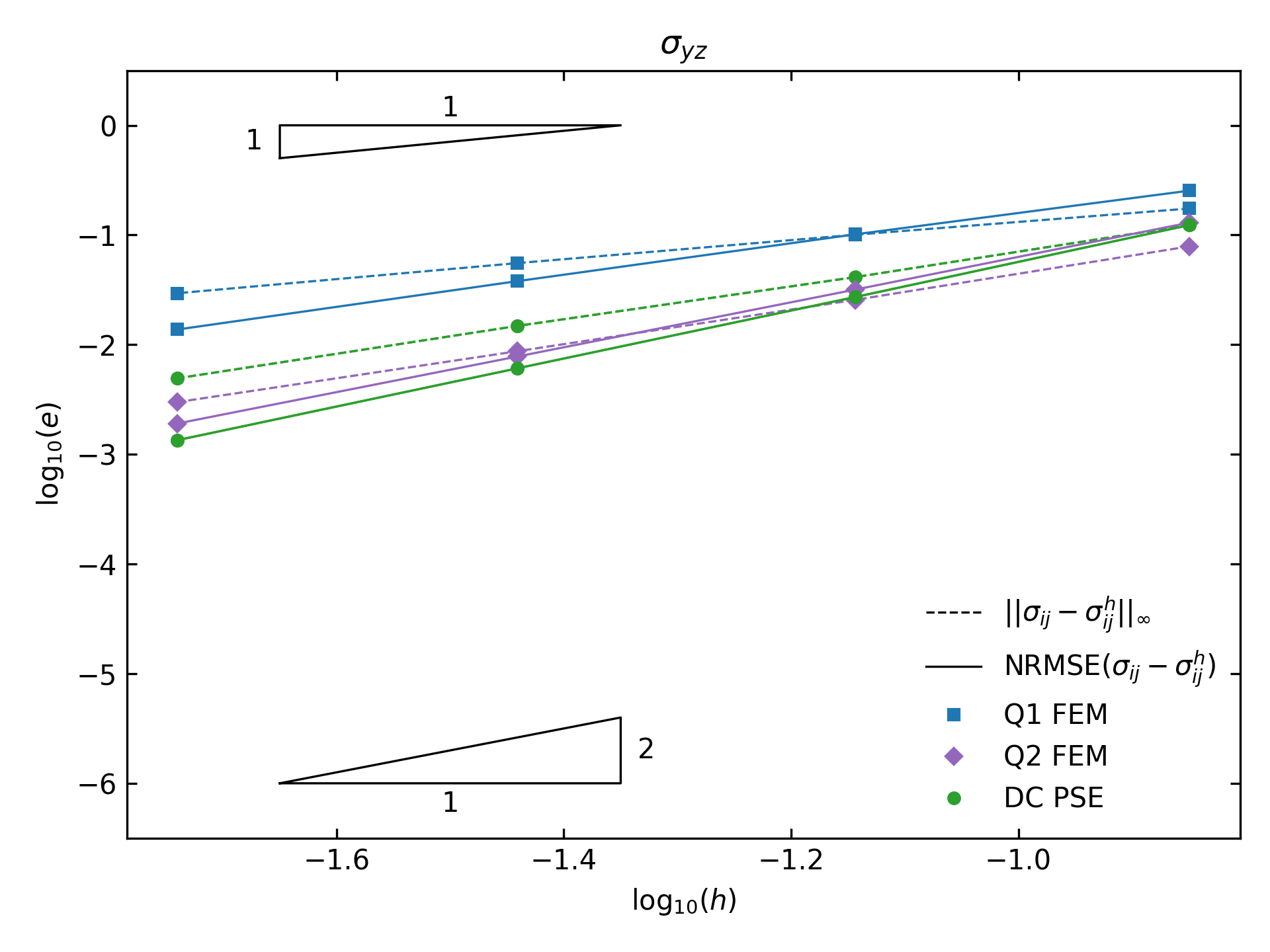}
    \caption{}
  \end{subfigure}
  \begin{subfigure}[b]{0.45\textwidth}
    \includegraphics[width=\textwidth]{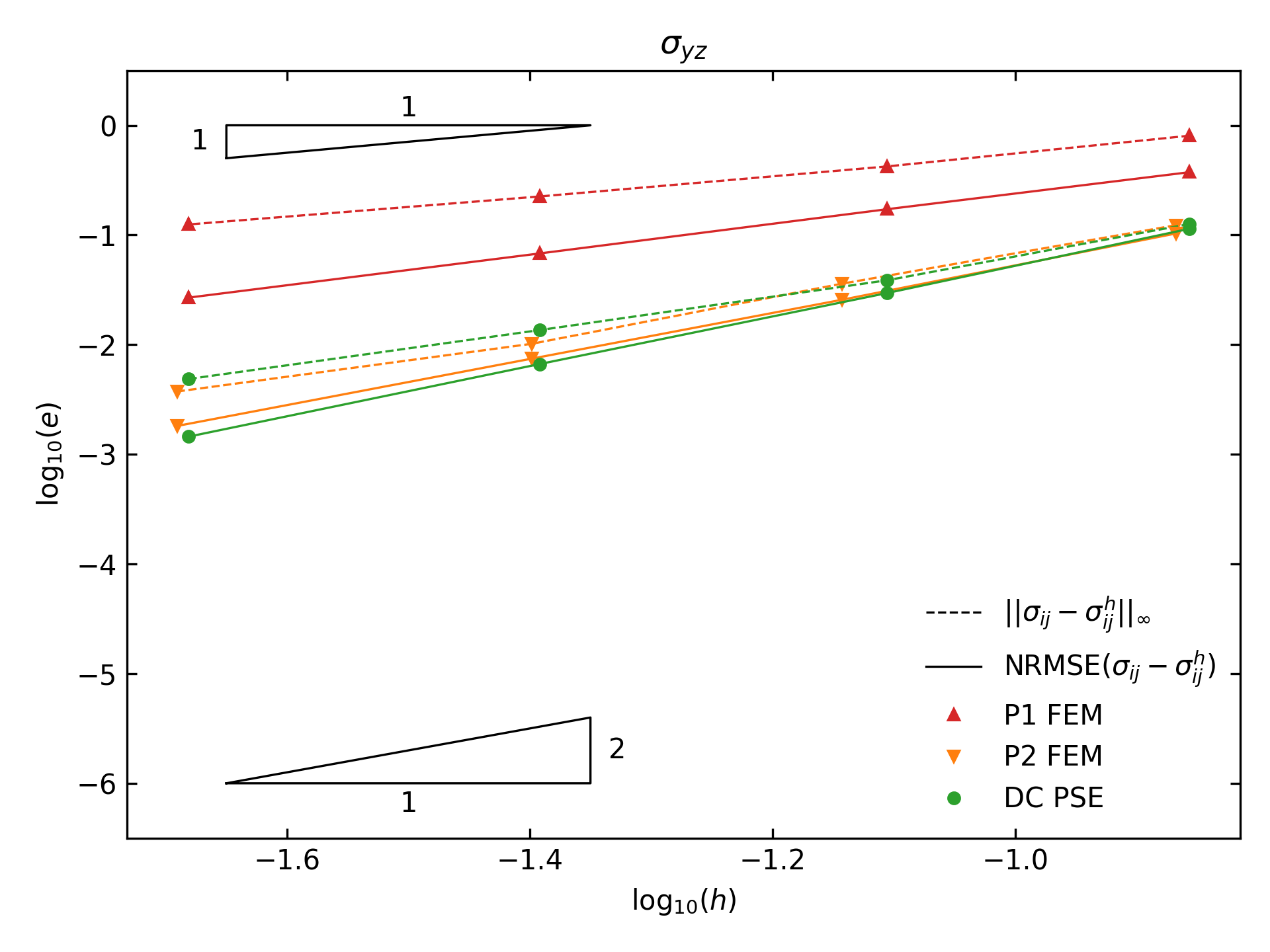}
    \caption{}
  \end{subfigure}  
  \\
  \begin{subfigure}[b]{0.45\textwidth}
    \includegraphics[width=\textwidth]{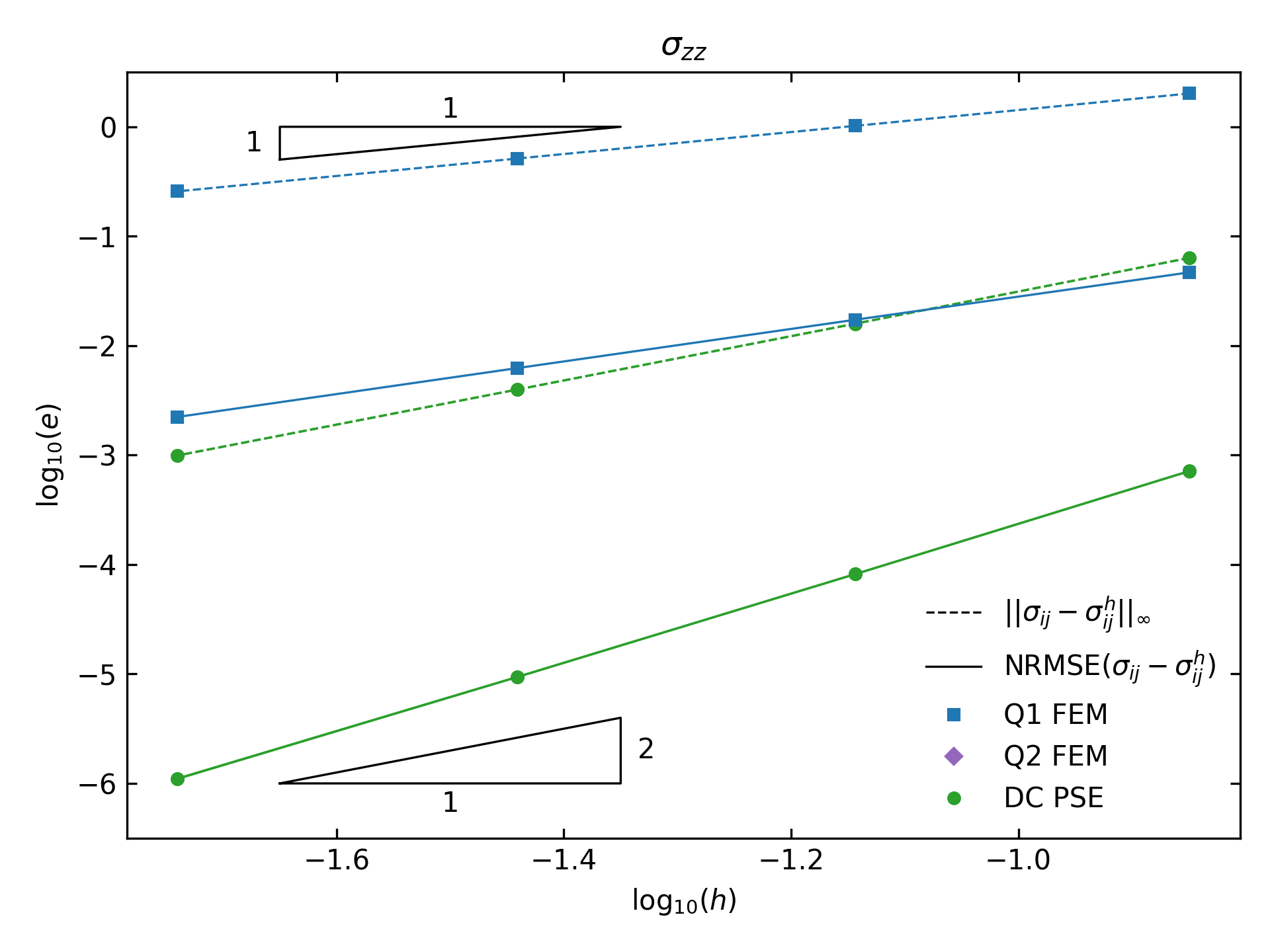}
    \caption{}
  \end{subfigure}
  \begin{subfigure}[b]{0.45\textwidth}
    \includegraphics[width=\textwidth]{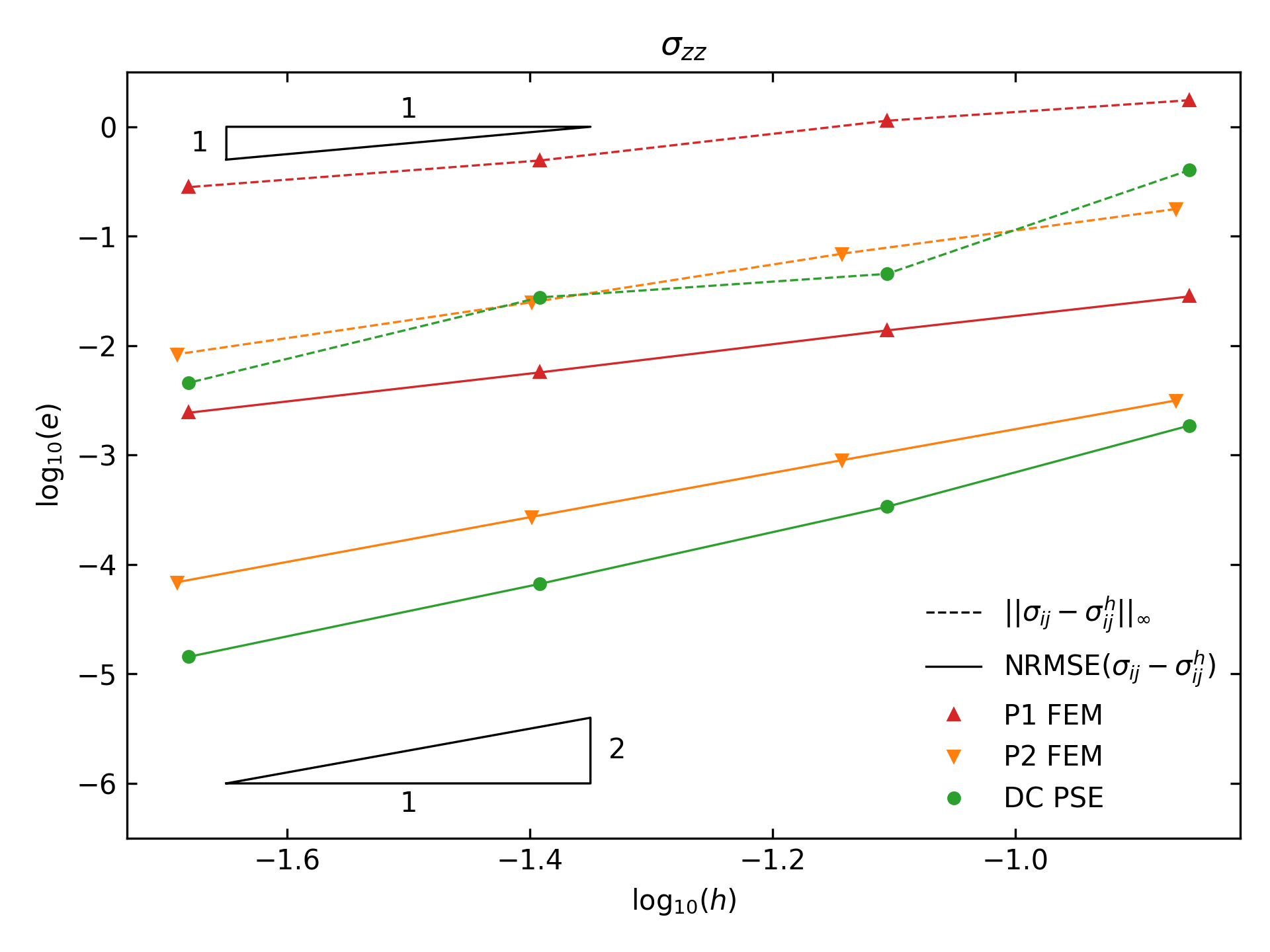}
    \caption{}
  \end{subfigure}  
  \caption{
    Convergence of stress components
    \(\sigma_{xz}\), \(\sigma_{yz}\) and \(\sigma_{zz}\)
    for the 3D cantilever benchmark problem
    using
    (a, c, e) structured hexahedral and
    (b, d, f) unstructured tetrahedral meshes.
    The normalized node spacing \(h\) for a 3D mesh with \(n\) nodes
    is defined as \(h = 1/(n^{1/3}-1)\).
  }
  \label{fig:cantilever-results}
\end{figure}

Fig.~\ref{fig:cantilever-meshes}
shows some representative meshes used for this problem.
Using the analytical solution for the displacements at the nodal points,
we computed the derivatives and derived quantities (strain and stress) numerically
using the FE method with nodal averaging
and the proposed DC PSE based meshfree recovery method.
Fig.~\ref{fig:cantilever-results} shows the errors in the nonzero stress components.
The proposed DC PSE recovery method performs well 
compared to the FE method,
with similar convergence rate and accuracy
as quadratic finite elements.
The $\sigma_{zz}$ stress component varies linearly
with respect to $y$ and $z$,
and is reproduced exactly by quadratic hexahedrons
(the error is less than $10^{-6}$ for all $h$ tested)
but not by linear finite elements or DC PSE.
For unstructured tetrahedral meshes
the DC PSE method appears to be slightly more accurate 
than the FE method with nodal averaging.
This suggests that DC PSE may be useful in practical applications
where complicated geometries cannot be meshed 
using a regular grid of hexahedral elements.

\FloatBarrier
\subsection{Stress in abdominal aortic aneurysm (AAA)}

An abdominal aortic aneurysm (AAA)
is an abnormal enlargement of the aorta,
which is the main blood vessel
that delivers blood to the body.
Rupture of an AAA can have fatal consequences
which means that patients must have periodic monitoring and, potentially,
surgical intervention to repair the aneurysm.
Stress analysis of aneurysms has been proposed as a risk assessment method for rupture
\citep{joldes_etal_2017_bioparr}.
The maximum stress often occurs on the inner or outer surface
of the aneurysm wall.
However, in the FE method,
the stress is typically computed at the integration points
which are located within the interior of the elements
and, moreover, the normal component of stress across element boundaries is discontinuous.
Stress recovery methods are required to obtain a smooth stress field
and accurate stresses on the surface nodes of the aneurysm wall.
In this example,
we demonstrate the application of our proposed approach
to the computation of stress in the arterial wall
given a pre-computed displacement field.
Here, we use the FEM to compute the displacement
and compare the stress computed using FEM
with that obtained using DC PSE.

The aneurysm geometry was extracted from a computed tomography (CT) scan
using an image-based segmentation process.
The open source mesh generator Gmsh \citep{geuzaine_remacle_2009_gmsh}
was used to create a mesh of the aneurysm wall
with 521,441 quadratic tetrahedral elements
and 928,682 nodes.
We used Altair HyperMesh
(\url{https://www.altair.com/hypermesh})
to check the element quality with the volumetric skew measure,
which is a volume ratio between an ideal shaped equilateral tetrahedron
and the tested actual tetrahedron.
The volumetric skew measure ranges between 0 and 1,
with 0 being an equilateral tetrahedron and 1 being a co-planar tetrahedron.
For the AAA mesh considered in this example,
all tetrahedra had volumetric skew less than 0.98.
To model the material behavior of the aneurysm wall
we used a linear elastic material with Young's modulus $E=60$ MPa 
and Poisson's ratio $\nu=0.45$
\citep{karimi_etal_2013_finite}.
We used Abaqus/Standard to compute the displacements and stresses using the finite element method
with quadratic tetrahedral displacement-based elements
(C3D10 element type in Abaqus).

\begin{figure}
  \centering
  \includegraphics[width=\textwidth]{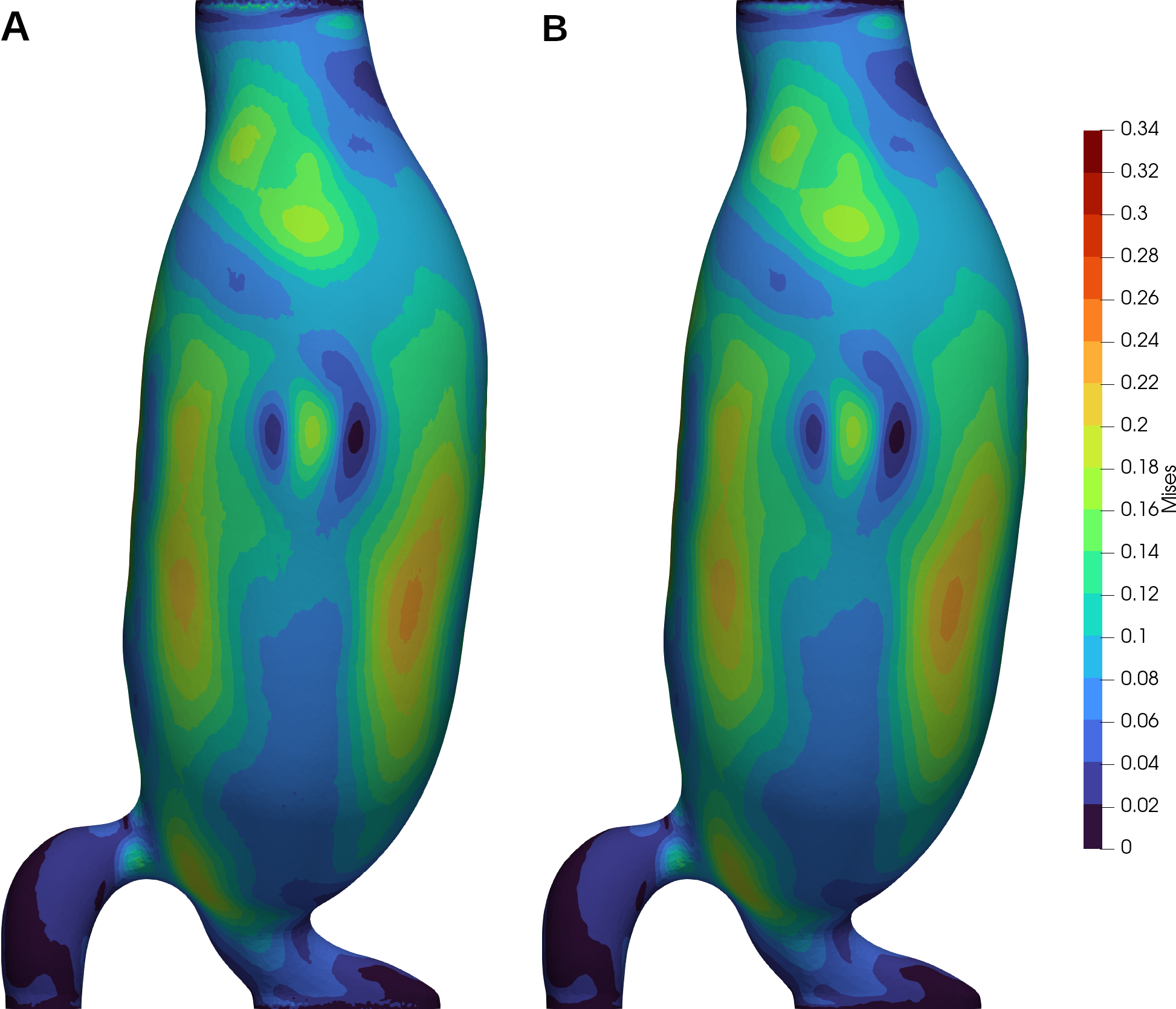}
  \caption{Distribution of von Mises stress
    on the outer wall of the abdominal aortic aneurysm (AAA).
    The displacement was computed using FEM
    with quadratic tetrahedral elements,
    and the stress was recovered using
    (a) FEM with nodal averaging and
    (b) meshfree DC PSE.
  }
  \label{fig:aaa-stress-distribution}
\end{figure}

Fig.~\ref{fig:aaa-stress-distribution}
shows the von Mises stress distribution computed using the FEM with nodal averaging compared to that computed using the proposed DC PSE meshfree approach.
The stress distributions computed using the two methods are very similar 
but upon close inspection it is apparent that the DC PSE solution
is smoother in some regions.
Table~\ref{tab:aaa-error}
shows the NRMSE between the stress components computed using FEM %
and those computed using the meshfree DC PSE based recovery method.
The difference between the stress components computed using the FEM and DC PSE is less than $1\%$.
Given that for the benchmark problems above
the stress computed using quadratic tetrahedral elements
with nodal averaging was slightly less accurate compared to DC PSE,
we may believe that the results for the stress in the AAA computed using DC PSE should also be more accurate than those computed using FEM with nodal averaging.

\begin{table}
    \caption{Normalized root mean square error (NRMSE)
    between the stress components computed using FEM %
    and those computed using the DC PSE based recovery method
    for the abdominal aortic aneurysm (AAA) example problem.}
    \label{tab:aaa-error}
    \centering
    \begin{tabular}{@{}lc@{}}
    \toprule
    Variable & NRMSE \\
    \midrule
    $\sigma_{xx}$ & 9.51$\times 10^{-3}$ \\
    $\sigma_{xy}$ & 2.38$\times 10^{-3}$ \\
    $\sigma_{xz}$ & 4.11$\times 10^{-3}$ \\
    $\sigma_{yy}$ & 9.53$\times 10^{-3}$ \\
    $\sigma_{yz}$ & 3.52$\times 10^{-3}$ \\
    $\sigma_{zz}$ & 7.78$\times 10^{-3}$ \\
    $\sigma_{\mathrm{Mises}}$ & 5.80$\times 10^{-3}$ \\
    \bottomrule
    \end{tabular}
\end{table}

\FloatBarrier
\section{Discussion}
\label{sec:discussion}

We presented a novel recovery method based on DC PSE derivative operators.
The DC PSE method is a truly meshfree method
that does not require nodal connectivity defined by elements
or the construction of element patches.
The recovery procedure requires only the nodal points and nodal values
of the primary variables to compute gradients and strains.
When the constitutive material model and its parameters are known
the stresses can also be computed.
This allows the method to be easily incorporated
into existing finite element software
or used as a standalone post-processing tool.

To assess the accuracy of the proposed approach
we compared its results with those obtained using the finite element method
for a number of benchmark problems,
before demonstrating the application of the method
to a practical problem with complicated geometry
by computing the stress within the wall of an aneurysm.
In most cases, the meshless approach
provides superior estimates of the gradient and other derived quantities
compared to the FE method with nodal averaging.
There are of course more sophisticated mesh-based methods
that can be used to obtain more accurate gradients than nodal averaging,
however, one of the salient advantages of the proposed method
is that the gradient can be computed with accuracy similar to quadratic finite elements
using a mesh of linear elements with a similar number of nodes.
This is because the nodal connectivity is not required
in our approach and the accuracy depends only on the density of the nodal distribution.
Therefore, if the solution (displacement)
can be computed with high accuracy using a dense mesh of linear elements,
its gradient (and derived quantities of strain and stress)
can be computed using the meshless DC PSE method
with similar accuracy as if a coarser mesh of quadratic elements
with a similar number of nodes
was used to compute the solution and its gradient.
This has significant practical importance in problems with complicated geometries
that are reconstructed from image data
with the boundaries triangulated using linear triangular elements.
In these cases it is often difficult or impractical
to generate a high quality mesh using higher order elements with curved edges
that accurately capture the geometry of the boundary.

The new recovery method based on the DC PSE formulation is general
and can be used for both linear and nonlinear problems.
Higher order derivatives than those used in the current study
can also be easily computed using the DC PSE method.
Although in this study we considered only the equations of linear elasticity,
the recovery method can be easily applied to other boundary value problems such as
nonlinear solid mechanics and the Navier--Stokes equations
by simply modifying the equations used to compute the derived quantities
given the gradient obtained using DC PSE.
An important quantity of interest in fluid flow simulations is the wall shear stress
\citep{gijsen_etal_2019_expert}.
In many practical applications the velocity field is obtained using piecewise linear (P1) elements
which makes accurate recovery of the gradient a difficult task
\citep{valen-sendstad_steinman_2014_mind}.
This task may be simplified by the adoption of the recovery procedure presented
and we intend to investigate the accuracy of this approach in our future studies.

The proposed recovery procedure may be applied as a recovery-based error estimator
\citep{ainsworth_oden_2000_posteriori}
whereby the discontinuous finite element solution at the element nodes
is compared to the continuous solution obtained using the DC PSE operators.
This error estimator could be used in adaptive mesh refinement algorithms
in both finite element and meshless methods.
The advantage of this approach over existing element-based methods
is that patches do not need to be created
thereby simplifying the implementation,
especially when using complex meshes consisting of various element shapes.

The meshfree gradient recovery approach that we have described in this paper
has many practical applications.
Its simplicity,
and the fact that it does not rely on predefined connectivity between neighboring nodes,
makes it an ideal post-processing tool
to be used as an ``add-on'' to existing simulation packages
for improving the accuracy of stress and strain computations.
We are now working on extending the approach described in this paper
to other equations with practical significance
with the aim of improving the estimation of derived quantities
from accurate solutions of the primary solution variable(s)
obtained using low order (typically linear) approximations
that are often employed in practice.

\section*{Acknowledgments}

This research was supported in part by
the Australian Government through
the Australian Research Council's Discovery Projects funding scheme
(project DP160100714).
The views expressed herein are those of the authors
and are not necessarily those of the Australian Research Council.

\newpage
\bibliography{manuscript}

\begin{thebibliography}{30}
\expandafter\ifx\csname natexlab\endcsname\relax\def\natexlab#1{#1}\fi
\providecommand{\url}[1]{\texttt{#1}}
\providecommand{\href}[2]{#2}
\providecommand{\path}[1]{#1}
\providecommand{\DOIprefix}{doi:}
\providecommand{\ArXivprefix}{arXiv:}
\providecommand{\URLprefix}{URL: }
\providecommand{\Pubmedprefix}{pmid:}
\providecommand{\doi}[1]{\href{http://dx.doi.org/#1}{\path{#1}}}
\providecommand{\Pubmed}[1]{\href{pmid:#1}{\path{#1}}}
\providecommand{\bibinfo}[2]{#2}
\ifx\xfnm\relax \def\xfnm[#1]{\unskip,\space#1}\fi
\bibitem[{Ahmed(2020)}]{ahmed_2020_comparative}
\bibinfo{author}{Ahmed, M.}, \bibinfo{year}{2020}.
\newblock \bibinfo{title}{A {{Comparative Study}} of {{Mesh-Free Radial Point
  Interpolation Method}} and {{Moving Least Squares Method-Based Error
  Estimation}} in {{Elastic Finite Element Analysis}}}.
\newblock \bibinfo{journal}{Arabian Journal for Science and Engineering}
  \bibinfo{volume}{45}, \bibinfo{pages}{3541--3557}.
\newblock \DOIprefix\doi{10.1007/s13369-019-04154-5}.
\bibitem[{Ahmed et~al.(2018)Ahmed, Singh and
  Desmukh}]{ahmed_etal_2018_interpolation}
\bibinfo{author}{Ahmed, M.}, \bibinfo{author}{Singh, D.},
  \bibinfo{author}{Desmukh, M.N.}, \bibinfo{year}{2018}.
\newblock \bibinfo{title}{Interpolation {{Type Stress Recovery Technique Based
  Error Estimator}} for {{Elasticity Problems}}}.
\newblock \bibinfo{journal}{Mechanics} \bibinfo{volume}{24},
  \bibinfo{pages}{672--679}.
\newblock \DOIprefix\doi{10.5755/j01.mech.24.5.19937}.
\bibitem[{Ainsworth and Oden(2000)}]{ainsworth_oden_2000_posteriori}
\bibinfo{author}{Ainsworth, M.}, \bibinfo{author}{Oden, J.T.},
  \bibinfo{year}{2000}.
\newblock \bibinfo{title}{A Posteriori Error Estimation in Finite Element
  Analysis}.
\newblock Pure and Applied Mathematics., \bibinfo{publisher}{{Wiley}},
  \bibinfo{address}{{New York}}.
\bibitem[{Barber(2010)}]{barber_2010_elasticity}
\bibinfo{author}{Barber, J.R.}, \bibinfo{year}{2010}.
\newblock \bibinfo{title}{Elasticity}. volume \bibinfo{volume}{172} of
  \textit{\bibinfo{series}{Solid {{Mechanics}} and {{Its Applications}}}}.
\newblock \bibinfo{publisher}{{Springer Netherlands}},
  \bibinfo{address}{{Dordrecht}}.
\newblock \DOIprefix\doi{10.1007/978-90-481-3809-8}.
\bibitem[{Bishop(2014)}]{bishop_2014_displacementbased}
\bibinfo{author}{Bishop, J.E.}, \bibinfo{year}{2014}.
\newblock \bibinfo{title}{A displacement-based finite element formulation for
  general polyhedra using harmonic shape functions}.
\newblock \bibinfo{journal}{International Journal for Numerical Methods in
  Engineering} \bibinfo{volume}{97}, \bibinfo{pages}{1--31}.
\newblock \DOIprefix\doi{10.1002/nme.4562}.
\bibitem[{Boroomand and
  Zienkiewicz(1997)}]{boroomand_zienkiewicz_1997_recovery}
\bibinfo{author}{Boroomand, B.}, \bibinfo{author}{Zienkiewicz, O.C.},
  \bibinfo{year}{1997}.
\newblock \bibinfo{title}{Recovery by {{Equilibrium}} in {{Patches}}
  ({{REP}})}.
\newblock \bibinfo{journal}{International Journal for Numerical Methods in
  Engineering} \bibinfo{volume}{40}, \bibinfo{pages}{137--164}.
\newblock
  \DOIprefix\doi{10.1002/(SICI)1097-0207(19970115)40:1<137::AID-NME57>3.0.CO;2-5}.
\bibitem[{Chen and Belytschko(2015)}]{chen_belytschko_2015_meshless}
\bibinfo{author}{Chen, J.S.}, \bibinfo{author}{Belytschko, T.},
  \bibinfo{year}{2015}.
\newblock \bibinfo{title}{Meshless and {{Meshfree Methods}}}, in:
  \bibinfo{editor}{Engquist, B.} (Ed.), \bibinfo{booktitle}{Encyclopedia of
  {{Applied}} and {{Computational Mathematics}}}. \bibinfo{publisher}{{Springer
  Berlin Heidelberg}}, \bibinfo{address}{{Berlin, Heidelberg}}, pp.
  \bibinfo{pages}{886--894}.
\newblock \DOIprefix\doi{10.1007/978-3-540-70529-1_531}.
\bibitem[{Choi et~al.(1991)Choi, Thorpe and Hanna}]{choi_etal_1991_image}
\bibinfo{author}{Choi, D.}, \bibinfo{author}{Thorpe, J.L.},
  \bibinfo{author}{Hanna, R.B.}, \bibinfo{year}{1991}.
\newblock \bibinfo{title}{Image analysis to measure strain in wood and paper}.
\newblock \bibinfo{journal}{Wood Science and Technology} \bibinfo{volume}{25},
  \bibinfo{pages}{251--262}.
\newblock \DOIprefix\doi{10.1007/BF00225465}.
\bibitem[{Fasshauer(2007)}]{fasshauer_2007_meshfree}
\bibinfo{author}{Fasshauer, G.E.}, \bibinfo{year}{2007}.
\newblock \bibinfo{title}{Meshfree {{Approximation Methods}} with {{MATLAB}}}.
\newblock \bibinfo{publisher}{{World Scientific}},
  \bibinfo{address}{{Singapore}}.
\bibitem[{Franke(1979)}]{franke_1979_critical}
\bibinfo{author}{Franke, R.}, \bibinfo{year}{1979}.
\newblock \bibinfo{title}{A {{Critical Comparison}} of {{Some Methods}} for
  {{Interpolation}} of {{Scattered Data}}}.
\newblock \bibinfo{type}{Technical {{Report}}}. {Monterey, California: Naval
  Postgraduate School.}
\newblock \URLprefix \url{https://calhoun.nps.edu/handle/10945/35052}.
\bibitem[{Geuzaine and Remacle(2009)}]{geuzaine_remacle_2009_gmsh}
\bibinfo{author}{Geuzaine, C.}, \bibinfo{author}{Remacle, J.F.},
  \bibinfo{year}{2009}.
\newblock \bibinfo{title}{Gmsh: {{A}} 3-{{D}} finite element mesh generator
  with built-in pre- and post-processing facilities}.
\newblock \bibinfo{journal}{International Journal for Numerical Methods in
  Engineering} \bibinfo{volume}{79}, \bibinfo{pages}{1309--1331}.
\newblock \DOIprefix\doi{10.1002/nme.2579}.
\bibitem[{Gijsen et~al.(2019)Gijsen, Katagiri, Barlis, Bourantas, Collet,
  Coskun, Daemen, Dijkstra, Edelman, Evans, {van der Heiden}, Hose, Koo, Krams,
  Marsden, Migliavacca, Onuma, Ooi, Poon, Samady, Stone, Takahashi, Tang,
  Thondapu, Tenekecioglu, Timmins, Torii, Wentzel and
  Serruys}]{gijsen_etal_2019_expert}
\bibinfo{author}{Gijsen, F.}, \bibinfo{author}{Katagiri, Y.},
  \bibinfo{author}{Barlis, P.}, \bibinfo{author}{Bourantas, C.},
  \bibinfo{author}{Collet, C.}, \bibinfo{author}{Coskun, U.},
  \bibinfo{author}{Daemen, J.}, \bibinfo{author}{Dijkstra, J.},
  \bibinfo{author}{Edelman, E.}, \bibinfo{author}{Evans, P.},
  \bibinfo{author}{{van der Heiden}, K.}, \bibinfo{author}{Hose, R.},
  \bibinfo{author}{Koo, B.K.}, \bibinfo{author}{Krams, R.},
  \bibinfo{author}{Marsden, A.}, \bibinfo{author}{Migliavacca, F.},
  \bibinfo{author}{Onuma, Y.}, \bibinfo{author}{Ooi, A.},
  \bibinfo{author}{Poon, E.}, \bibinfo{author}{Samady, H.},
  \bibinfo{author}{Stone, P.}, \bibinfo{author}{Takahashi, K.},
  \bibinfo{author}{Tang, D.}, \bibinfo{author}{Thondapu, V.},
  \bibinfo{author}{Tenekecioglu, E.}, \bibinfo{author}{Timmins, L.},
  \bibinfo{author}{Torii, R.}, \bibinfo{author}{Wentzel, J.},
  \bibinfo{author}{Serruys, P.}, \bibinfo{year}{2019}.
\newblock \bibinfo{title}{Expert recommendations on the assessment of wall
  shear stress in human coronary arteries: Existing methodologies, technical
  considerations, and clinical applications}.
\newblock \bibinfo{journal}{European Heart Journal} \bibinfo{volume}{40},
  \bibinfo{pages}{3421--3433}.
\newblock \DOIprefix\doi{10.1093/eurheartj/ehz551}.
\bibitem[{Joldes et~al.(2017)Joldes, Miller, Wittek, Forsythe, Newby and
  Doyle}]{joldes_etal_2017_bioparr}
\bibinfo{author}{Joldes, G.R.}, \bibinfo{author}{Miller, K.},
  \bibinfo{author}{Wittek, A.}, \bibinfo{author}{Forsythe, R.O.},
  \bibinfo{author}{Newby, D.E.}, \bibinfo{author}{Doyle, B.J.},
  \bibinfo{year}{2017}.
\newblock \bibinfo{title}{{{BioPARR}}: {{A}} software system for estimating the
  rupture potential index for abdominal aortic aneurysms}.
\newblock \bibinfo{journal}{Scientific Reports} \bibinfo{volume}{7},
  \bibinfo{pages}{4641}.
\newblock \DOIprefix\doi{10.1038/s41598-017-04699-1}.
\bibitem[{Karimi et~al.(2013)Karimi, Navidbakhsh, Faghihi, Shojaei and
  Hassani}]{karimi_etal_2013_finite}
\bibinfo{author}{Karimi, A.}, \bibinfo{author}{Navidbakhsh, M.},
  \bibinfo{author}{Faghihi, S.}, \bibinfo{author}{Shojaei, A.},
  \bibinfo{author}{Hassani, K.}, \bibinfo{year}{2013}.
\newblock \bibinfo{title}{A finite element investigation on plaque
  vulnerability in realistic healthy and atherosclerotic human coronary
  arteries}.
\newblock \bibinfo{journal}{Proceedings of the Institution of Mechanical
  Engineers, Part H: Journal of Engineering in Medicine} \bibinfo{volume}{227},
  \bibinfo{pages}{148--161}.
\newblock \DOIprefix\doi{10.1177/0954411912461239}.
\bibitem[{Kelly et~al.(1983)Kelly, Gago, Zienkiewicz and
  Babuska}]{kelly_etal_1983_posteriori}
\bibinfo{author}{Kelly, D.W.}, \bibinfo{author}{Gago, J.P.D.S.R.},
  \bibinfo{author}{Zienkiewicz, O.C.}, \bibinfo{author}{Babuska, I.},
  \bibinfo{year}{1983}.
\newblock \bibinfo{title}{A posteriori error analysis and adaptive processes in
  the finite element method: {{Part I}}\textemdash error analysis}.
\newblock \bibinfo{journal}{International Journal for Numerical Methods in
  Engineering} \bibinfo{volume}{19}, \bibinfo{pages}{1593--1619}.
\newblock \DOIprefix\doi{10.1002/nme.1620191103}.
\bibitem[{Kirsch(1898)}]{kirsch_1898_theorie}
\bibinfo{author}{Kirsch, G.}, \bibinfo{year}{1898}.
\newblock \bibinfo{title}{{Die Theorie der Elastizit\"at und die Bed\"urfnisse
  der Festigkeitslehre}}.
\newblock \bibinfo{journal}{Zeitschrift des Vereines Deutscher Ingenieure}
  \bibinfo{volume}{42}, \bibinfo{pages}{797--807}.
\bibitem[{Lee and Zhou(2004)}]{lee_zhou_2004_error}
\bibinfo{author}{Lee, C.K.}, \bibinfo{author}{Zhou, C.E.},
  \bibinfo{year}{2004}.
\newblock \bibinfo{title}{On error estimation and adaptive refinement for
  element free {{Galerkin}} method: {{Part I}}: Stress recovery and a
  posteriori error estimation}.
\newblock \bibinfo{journal}{Computers \& Structures} \bibinfo{volume}{82},
  \bibinfo{pages}{413--428}.
\newblock \DOIprefix\doi{10.1016/j.compstruc.2003.10.018}.
\bibitem[{{L{\'o}pez-Linares} et~al.(2019){L{\'o}pez-Linares}, Garc{\'i}a,
  Garc{\'i}a, Cortes, Piella, Mac{\'i}a, Noailly and
  Gonz{\'a}lez~Ballester}]{lopez-linares_etal_2019_imagebased}
\bibinfo{author}{{L{\'o}pez-Linares}, K.}, \bibinfo{author}{Garc{\'i}a, I.},
  \bibinfo{author}{Garc{\'i}a, A.}, \bibinfo{author}{Cortes, C.},
  \bibinfo{author}{Piella, G.}, \bibinfo{author}{Mac{\'i}a, I.},
  \bibinfo{author}{Noailly, J.}, \bibinfo{author}{Gonz{\'a}lez~Ballester,
  M.A.}, \bibinfo{year}{2019}.
\newblock \bibinfo{title}{Image-{{Based 3D Characterization}} of {{Abdominal
  Aortic Aneurysm Deformation After Endovascular Aneurysm Repair}}}.
\newblock \bibinfo{journal}{Frontiers in Bioengineering and Biotechnology}
  \bibinfo{volume}{7}.
\newblock \DOIprefix\doi{10.3389/fbioe.2019.00267}.
\bibitem[{Oden and Brauchli(1971)}]{oden_brauchli_1971_calculation}
\bibinfo{author}{Oden, J.T.}, \bibinfo{author}{Brauchli, H.J.},
  \bibinfo{year}{1971}.
\newblock \bibinfo{title}{On the calculation of consistent stress distributions
  in finite element approximations}.
\newblock \bibinfo{journal}{International Journal for Numerical Methods in
  Engineering} \bibinfo{volume}{3}, \bibinfo{pages}{317--325}.
\newblock \DOIprefix\doi{10.1002/nme.1620030303}.
\bibitem[{Reboux et~al.(2012)Reboux, Schrader and
  Sbalzarini}]{reboux_etal_2012_selforganizing}
\bibinfo{author}{Reboux, S.}, \bibinfo{author}{Schrader, B.},
  \bibinfo{author}{Sbalzarini, I.F.}, \bibinfo{year}{2012}.
\newblock \bibinfo{title}{A self-organizing {{Lagrangian}} particle method for
  adaptive-resolution advection\textendash diffusion simulations}.
\newblock \bibinfo{journal}{Journal of Computational Physics}
  \bibinfo{volume}{231}, \bibinfo{pages}{3623--3646}.
\newblock \DOIprefix\doi{10.1016/j.jcp.2012.01.026}.
\bibitem[{Schrader(2011)}]{schrader_2011_discretizationcorrected}
\bibinfo{author}{Schrader, B.}, \bibinfo{year}{2011}.
\newblock \bibinfo{title}{Discretization-{{Corrected PSE Operators}} for
  {{Adaptive Multiresolution Particle Methods}}}.
\newblock \bibinfo{type}{Doctoral {{Thesis}}}. ETH Zurich.
\newblock \DOIprefix\doi{10.3929/ethz-a-006425176}.
\bibitem[{Schrader et~al.(2010)Schrader, Reboux and
  Sbalzarini}]{schrader_etal_2010_discretization}
\bibinfo{author}{Schrader, B.}, \bibinfo{author}{Reboux, S.},
  \bibinfo{author}{Sbalzarini, I.F.}, \bibinfo{year}{2010}.
\newblock \bibinfo{title}{Discretization correction of general integral {{PSE
  Operators}} for particle methods}.
\newblock \bibinfo{journal}{Journal of Computational Physics}
  \bibinfo{volume}{229}, \bibinfo{pages}{4159--4182}.
\newblock \DOIprefix\doi{10.1016/j.jcp.2010.02.004}.
\bibitem[{Schrader et~al.(2012)Schrader, Reboux and
  Sbalzarini}]{schrader_etal_2012_choosing}
\bibinfo{author}{Schrader, B.}, \bibinfo{author}{Reboux, S.},
  \bibinfo{author}{Sbalzarini, I.F.}, \bibinfo{year}{2012}.
\newblock \bibinfo{title}{Choosing the {{Best Kernel}}: {{Performance Models}}
  for {{Diffusion Operators}} in {{Particle Methods}}}.
\newblock \bibinfo{journal}{SIAM Journal on Scientific Computing}
  \bibinfo{volume}{34}, \bibinfo{pages}{A1607--A1634}.
\newblock \DOIprefix\doi{10.1137/110835815}.
\bibitem[{Timoshenko and Goodier(1951)}]{timoshenko_goodier_1951_theory}
\bibinfo{author}{Timoshenko, S.}, \bibinfo{author}{Goodier, J.N.},
  \bibinfo{year}{1951}.
\newblock \bibinfo{title}{Theory of {{Elasticity}}}.
\newblock \bibinfo{edition}{Second} ed., \bibinfo{publisher}{{McGraw-Hill}},
  \bibinfo{address}{{New York}}.
\bibitem[{{Valen-Sendstad} and
  Steinman(2014)}]{valen-sendstad_steinman_2014_mind}
\bibinfo{author}{{Valen-Sendstad}, K.}, \bibinfo{author}{Steinman, D.A.},
  \bibinfo{year}{2014}.
\newblock \bibinfo{title}{Mind the {{Gap}}: {{Impact}} of {{Computational Fluid
  Dynamics Solution Strategy}} on {{Prediction}} of {{Intracranial Aneurysm
  Hemodynamics}} and {{Rupture Status Indicators}}}.
\newblock \bibinfo{journal}{American Journal of Neuroradiology}
  \bibinfo{volume}{35}, \bibinfo{pages}{536--543}.
\newblock \DOIprefix\doi{10.3174/ajnr.A3793}.
\bibitem[{Zhang and Naga(2005)}]{zhang_naga_2005_new}
\bibinfo{author}{Zhang, Z.}, \bibinfo{author}{Naga, A.}, \bibinfo{year}{2005}.
\newblock \bibinfo{title}{A {{New Finite Element Gradient Recovery Method}}:
  {{Superconvergence Property}}}.
\newblock \bibinfo{journal}{SIAM Journal on Scientific Computing}
  \bibinfo{volume}{26}, \bibinfo{pages}{1192--1213}.
\newblock \DOIprefix\doi{10.1137/S1064827503402837}.
\bibitem[{Zienkiewicz et~al.(2013)Zienkiewicz, Taylor and
  Fox}]{zienkiewicz_etal_2013_finite_solid}
\bibinfo{author}{Zienkiewicz, O.C.}, \bibinfo{author}{Taylor, R.L.},
  \bibinfo{author}{Fox, D.D.}, \bibinfo{year}{2013}.
\newblock \bibinfo{title}{The {{Finite Element Method}} for {{Solid}} and
  {{Structural Mechanics}}}.
\newblock \bibinfo{edition}{Seventh} ed.,
  \bibinfo{publisher}{{Butterworth-Heinemann}}.
\bibitem[{Zienkiewicz and Zhu(1987)}]{zienkiewicz_zhu_1987_simple}
\bibinfo{author}{Zienkiewicz, O.C.}, \bibinfo{author}{Zhu, J.Z.},
  \bibinfo{year}{1987}.
\newblock \bibinfo{title}{A simple error estimator and adaptive procedure for
  practical engineerng analysis}.
\newblock \bibinfo{journal}{International Journal for Numerical Methods in
  Engineering} \bibinfo{volume}{24}, \bibinfo{pages}{337--357}.
\newblock \DOIprefix\doi{10.1002/nme.1620240206}.
\bibitem[{Zienkiewicz and
  Zhu(1992)}]{zienkiewicz_zhu_1992_superconvergent_part1_recovery}
\bibinfo{author}{Zienkiewicz, O.C.}, \bibinfo{author}{Zhu, J.Z.},
  \bibinfo{year}{1992}.
\newblock \bibinfo{title}{The superconvergent patch recovery and a posteriori
  error estimates. {{Part}} 1: {{The}} recovery technique}.
\newblock \bibinfo{journal}{International Journal for Numerical Methods in
  Engineering} \bibinfo{volume}{33}, \bibinfo{pages}{1331--1364}.
\newblock \DOIprefix\doi{10.1002/nme.1620330702}.
\bibitem[{Zienkiewicz and Zhu(1995)}]{zienkiewicz_zhu_1995_superconvergence}
\bibinfo{author}{Zienkiewicz, O.C.}, \bibinfo{author}{Zhu, J.Z.},
  \bibinfo{year}{1995}.
\newblock \bibinfo{title}{Superconvergence and the superconvergent patch
  recovery}.
\newblock \bibinfo{journal}{Finite Elements in Analysis and Design}
  \bibinfo{volume}{19}, \bibinfo{pages}{11--23}.
\newblock \DOIprefix\doi{10.1016/0168-874X(94)00054-J}.

\end{thebibliography}

\nolinenumbers

\newpage
\appendix

\end{document}